\documentclass[letterpaper, 10 pt, conference]{ieeeconf}  

\IEEEoverridecommandlockouts                              

\usepackage{amsmath}
\usepackage{amssymb}
\usepackage{amsthm}
\usepackage[bookmarksopen,bookmarksnumbered,colorlinks=true, linkcolor=black,citecolor = blue, urlcolor = blue, filecolor=blue, pagebackref=false, hypertexnames=false, plainpages=false, pdfpagelabels]{hyperref}
\usepackage{subcaption}
\usepackage{graphicx}
\usepackage{wrapfig}
\usepackage{cite}
\usepackage{multirow}

\usepackage{tikz} 
\usepackage{tkz-graph}
\usepackage{tkz-berge}
\usetikzlibrary{backgrounds,fit,shapes,snakes,arrows,shapes.geometric,positioning}
\usetikzlibrary{intersections,patterns,shapes.misc}
\usetikzlibrary{decorations.pathmorphing}

\tikzstyle{block} = [rectangle, rounded corners, minimum width=3cm, minimum height=1cm,text centered, draw=black, fill=red!30]
\tikzstyle{new} = [rectangle, rounded corners, minimum width=1cm, minimum
height=1cm,text centered, draw=black, fill=blue!10!white, dashed]
\tikzstyle{arrow} = [thick,->,>=stealth]
\usetikzlibrary{calc, quotes}

\tikzstyle{fblock} = [rectangle, draw, fill=gray!20, 
text width=8em, text centered, rounded corners, minimum height=4em, minimum width = 8em]
\tikzstyle{line} = [draw, -latex']

\usepackage{algorithmicx}
\usepackage{algorithm}
\usepackage[]{algpseudocode}
\algnewcommand\algorithmicinput{\textbf{Input:}}
\algnewcommand\INPUT{\item[\algorithmicinput]}

\usepackage{fontawesome}

\usepackage{soul}
\usepackage{comment}

\newtheorem{theorem}{Theorem}
\newtheorem{lemma}{Lemma}

\newtheorem{problem}{Problem}
\newtheorem{remark}{Remark}

\newtheorem{definition}{Definition}
\newtheorem{assumption}{Assumption}

\DeclareMathOperator{\tr}{tr}

\DeclareMathOperator{\Stiefel}{St}

\DeclareMathOperator{\rgrad}{grad}

\DeclareMathOperator{\proj}{Proj}
\DeclareMathOperator{\Retr}{Retr}
\DeclareMathOperator{\Precon}{Precon}

\newcommand{\norm}[1]{\left\lVert#1\right\rVert}

\DeclareMathOperator{\Orthogonal}{O}
\DeclareMathOperator{\SOd}{SO}
\DeclareMathOperator{\SE}{SE}

\DeclareMathOperator{\Langevin}{Langevin}

\newcommand{\Nbr}{\text{N}}
\newcommand{\Mcal}{\mathcal{M}}
\newcommand{\Ncal}{\mathcal{N}}

\newcommand{\Ecal}{\mathcal{E}}

\newcommand{\Rcal}{\mathcal{R}}
\newcommand{\fhat}{\hat{f}}
\newcommand{\ik}{i_k}
\newcommand{\jk}{j_k}

\newcommand{\ek}{e_k}

\newcommand{\ALGORITHM}{\textsc{ASAPP}}

\newcommand \remove[1]      {{\color{red}}}

\newcommand \blue[1]        {{\color{black}#1}}

\title{\LARGE \bf
Asynchronous and Parallel Distributed Pose Graph Optimization
}

\author{Yulun Tian$^{1}$, Alec Koppel$^{2}$, Amrit Singh Bedi$^{2}$, and Jonathan P. How$^{1}$
\thanks{Supported in part by ARL
DCIST under Cooperative Agreement Number W911NF-17-2-0181
and by NASA Convergent Aeronautics Solutions project Design
Environment for Novel Vertical Lift Vehicles (DELIVER).}
\thanks{$^{1}$Y. Tian and J. P. How are with the Department of Aeronautics and Astronautics,  Massachusetts Institute of Technology, 77 Massachusetts Ave, Cambridge, MA, {\tt\small \{yulun,jhow\}@mit.edu}}%
\thanks{$^{2}$A. Koppel and A. S. Bedi are with the U.S. Army Research Laboratory, Adelphi, MD 20783
{\tt\small alec.e.koppel.civ@mail.mil, amrit0714@gmail.com}}%
}

\begin{document}

\maketitle
\thispagestyle{empty}
\pagestyle{empty}

\begin{abstract}

We present Asynchronous Stochastic Parallel Pose Graph Optimization ($\ALGORITHM$), the first \emph{asynchronous} algorithm for distributed pose graph optimization (PGO) in multi-robot simultaneous localization and mapping. 
By enabling robots to optimize their local trajectory estimates without synchronization, 
$\ALGORITHM$ offers resiliency against communication delays and alleviates the need to wait for stragglers in the network. 
\blue{Furthermore, $\ALGORITHM$ can be applied on the rank-restricted relaxations of PGO,
a crucial class of non-convex Riemannian optimization problems that underlies recent breakthroughs on globally optimal PGO.}
Under bounded delay, we establish the {global} first-order convergence of $\ALGORITHM$ using a sufficiently small stepsize.
The derived stepsize depends on the worst-case delay and inherent problem sparsity, and furthermore matches known result for synchronous algorithms when there is no delay. 
Numerical evaluations on simulated and real-world datasets demonstrate 
\blue{favorable performance compared to state-of-the-art synchronous approach, and show $\ALGORITHM$'s resilience against a wide range of delays in practice.}



\end{abstract}



\section{Introduction}
Multi-robot simultaneous localization and mapping (SLAM) is a fundamental capability for many real-world robotic applications.
\emph{Pose graph optimization} (PGO) is the backbone of state of the art approaches to multi-robot SLAM, which fuses individual trajectories together and endows participating robots with a common spatial understanding of the environment.
Many approaches to multi-robot PGO require the centralized processing of observations at a base station, which is communication intensive and vulnerable to single point of failure.
In contrast, decentralized approaches are favorable as they effectively mitigate communication, privacy, and vulnerability concerns associated with centralization.

Recent works on distributed PGO have achieved important progress; see e.g., \cite{Choudhary2017Distributed,Tian2019RBCD} and the references therein. 
However, to the best of our knowledge, existing distributed algorithms are inherently \emph{synchronous}, which necessitates that robots, for instance, pass messages over the network or wait at predetermined points, in order to ensure up-to-date information sharing during distributed optimization. Doing so may incur considerable communication overhead and increase the complexity of implementation. On the other hand, simply dropping synchronization in the execution of synchronous algorithms may slow down convergence or even cause divergence, both in theory and practice.

In this work, we overcome the aforementioned challenge by proposing $\ALGORITHM$ (\textbf{A}synchronous \textbf{S}toch\textbf{A}stic \textbf{P}arallel \textbf{P}ose Graph Optimization), the first \emph{asynchronous} and \emph{provably convergent} algorithm for distributed PGO. 
We take inspiration from existing parallel and asynchronous algorithms \cite{bertsekas1989parallel,Agarwal2011NIPS,Niu2011NIPS,Liu2015SIAM,Lian2015NIPS}, and adapt these ideas to solve the \emph{non-convex} Riemannian optimization problem underlying PGO.
In $\ALGORITHM$, each robot executes its local optimization loop at a high rate, without waiting for updates from others over the network. 
This makes $\ALGORITHM$ easier to implement in practice and flexible against communication delay.
\blue{
In addition, recent breakthroughs in centralized PGO, starting with SE-Sync~\cite{RosenSESyncIJRR19}, show how one can obtain globally optimal solutions to PGO by solving a hierarchy of rank-restricted relaxations.
In this work, we leverage the important insights provided by SE-Sync and subsequent works \cite{Briales2017RAL,Tian2019RBCD}, and show that $\ALGORITHM$ can solve both PGO and its rank-restricted relaxations.
}

\blue{
This work focuses on solving the distributed SLAM back-end (i.e., pose graph optimization).
Our solver can be combined with a separate front-end module into a full SLAM solution, similar to what is done in state-of-the-art multi-robot SLAM systems, e.g.,  \cite{Cieslewski2018ICRA,Lajoie2019RAL}.
}

{\bf \noindent Contributions}
Since asynchronous algorithms allow communication delays to be substantial and unpredictable,
it is usually unclear under what conditions they converge in practice. 
In this work, we provide a rigorous answer to this question and establish the first known convergence result for asynchronous algorithms on the \emph{non-convex} PGO problem.
In particular, we show that as long as the worst-case delay is not arbitrarily large, $\ALGORITHM$ 
\blue{with a sufficiently small stepsize always converges to first-order critical points when solving PGO and its rank-restricted relaxations, with \emph{global} sublinear convergence rate.}
The derived stepsize depends on the maximum delay and inherent problem sparsity, 
and furthermore reduces to the well known constant of $1/L$ (where $L$ is the Lipschitz constant) for synchronous algorithms when there is no delay.
Numerical evaluations on simulated and real-world datasets demonstrate that $\ALGORITHM$
\blue{compares favorably against state-of-the-art synchronous methods,}
and furthermore is resilient against a wide range of communication delays. 
Both results show the practical value of the proposed algorithm in a realistic distributed setting. 



\subsection*{Preliminaries on Riemannian Optimization}
\label{sec:preliminary}

This work relies heavily on the first-order geometry of Riemannian manifolds.
The reader is referred to \cite{absil2009optimization} for a rigorous treatment of this subject.
In SLAM, examples of matrix manifolds that frequently appear include the orthogonal group $\Orthogonal(d)$, special orthogonal group $\SOd(d)$, and the special Euclidean group $\SE(d)$.
In this work, we use $\Mcal \subseteq \Ecal$ to denote a general matrix submanifold, where $\Ecal$ is the so-called ambient space (in this work, $\Ecal$ is always the Euclidean space). 
Each point $x \in \Mcal$ on the manifold has an associated tangent space $T_x \Mcal$.
Informally, $T_x \Mcal$ contains all possible directions of change at $x$ while staying on $\Mcal$.
As $T_x \Mcal$ is a vector space, we also endow it with the standard Frobenius inner product, i.e., for two tangent vectors $\eta_1, \eta_2 \in T_x \Mcal$, $\langle \eta_1, \eta_2 \rangle \triangleq \tr(\eta_1^\top \eta_2)$. The inner product induces a norm $\norm{\eta} \triangleq \sqrt{\langle \eta, \eta \rangle}$.
Finally, a tangent vector can be mapped back to the manifold through a retraction $\Retr_x: T_x \Mcal \to \Mcal$, which is a smooth mapping that preserves the first-order structure of the manifold \cite{absil2009optimization}.

Riemannian optimization considers minimizing a function $f: \Mcal \to \mathbb{R}$ on the manifold. 
First-order Riemannian optimization algorithms, including the one proposed in this work, often use the Riemannian gradient $\rgrad f(x) \in T_x \Mcal$, which corresponds to the direction of steepest ascent in the tangent space. 
For matrix submanifolds, the Riemannian gradient is obtained by an orthogonal projection of the usual Euclidean gradient $\nabla f(x)$ onto the tangent space, i.e., $\rgrad f(x) = \proj_{T_x \Mcal} \nabla f(x)$ \cite{absil2009optimization}. 
We call $x^\star \in \Mcal$ a first-order critical point if $\rgrad f(x^\star) = 0$.


\section{Related Work}
\label{sec:related_work}

\subsection{Distributed and Parallel PGO}
In pursuit of decentralized \emph{asynchronous} algorithms, we note that synchronized decentralized PGO has been well-studied.
\blue{
Tron et al. \cite{Tron2009CDC,Tron2014} propose distributed Riemannian gradient descent on a set of reshape cost functions based on the geodesic distance, under which the method provably converges.
}
A similar gradient-based method with line-search has been proposed \cite{Knuth2012ICRA}.
\blue{
Cunningham et al.~\cite{CunninghamDDFSAM,Cunningham2012ICRA,CunninghamDDFSAM2} propose DDF-SAM, a decentralized system for feature-based SLAM where robots exchange Gaussian marginals over commonly observed variables. The use of local cache in \cite{CunninghamDDFSAM,Cunningham2012ICRA,CunninghamDDFSAM2} further allows asynchronous communication, although convergence is not discussed in this setting. }
Choudhary et al. \cite{Choudhary2015IROS} propose the alternating direction method of multipliers (ADMM) as a decentralized method to solve PGO.  
However, convergence of ADMM is not established due to the non-convex nature of the optimization problem. 
More recently, Choudhary et al. \cite{Choudhary2017Distributed} propose a two-stage approach
where each stage uses distributed Gauss-Seidel \cite{bertsekas1989parallel}
to solve a relaxed or linearized PGO problem.
The two-stage approach \cite{Choudhary2017Distributed} is further combined with outlier rejection schemes in \cite{Lajoie2019RAL}.
In our recent work \cite{Tian2019RBCD}, we avoid explicit linearization by directly optimizing PGO and its rank-restricted relaxations \cite{Briales2017RAL}. 
The proposed solver performs distributed block-coordinate descent over the product of Riemannian manifolds, and provably converges to first-order critical points with global sublinear rate. 
In a separate line of research, Fan and Murphey~\cite{Fan2019ISRR} propose an empirically accelerated PGO solver suitable for distributed optimization based on generalized proximal methods.

\subsection{Asynchronous Parallel Optimization}

The aforementioned works are promising, but critically rely on \emph{synchronization} for convergence, which limits their practical value for networked autonomous systems. However, within the broader optimization literature, there is a plethora of works on parallel and asynchronous optimization, partially motivated by popular applications in large-scale machine learning and deep learning. Study of asynchronous gradient-based algorithms began with the seminal work of Bertsekas and Tsitsilis~\cite{bertsekas1989parallel}, and have led to the recent development of asynchronous randomized block coordinate and stochastic gradient algorithms, see \cite{Agarwal2011NIPS,Niu2011NIPS,Liu2015SIAM,Liu2015JMLR,Lian2015NIPS,Lian2018ICML,Cannelli2019} and references therein. We are especially interested in asynchronous parallel schemes for non-convex optimization, which have been studied in \cite{Lian2015NIPS,Cannelli2019}. In this work, we generalize these approaches to the setting where the feasible set is the product of non-convex matrix manifolds, motivated by PGO.
\remove{
Our model of asynchrony is comparable to \cite{Lian2018ICML}: workers exchange local parameters asynchronously during optimization. 
However, unlike \cite{Lian2018ICML}, we obviate the need for local averaging to achieve consensus, as each robot is only responsible for updating its own trajectory.  }

\section{Problem Formulation}
\label{sec:problem_formulation}
\label{sec:pose_graph}


\begin{figure}[t]
	\centering
	\begin{subfigure}[t]{0.24\textwidth}
		\centering
		\begin{tikzpicture}[scale=1.1]
		\tikzstyle{vertex}=[circle,fill,scale=0.4,draw]
		\tikzstyle{special vertex}=[circle,fill=gray!30,scale=0.4,draw]
		\tikzstyle{square vertex}=[rectangle,fill,scale=0.5,draw]
		\tikzstyle{diamond vertex}=[regular polygon,regular polygon
		sides=3,rotate=45,fill,scale=0.3,draw]
		\node[special vertex] at (1,   2.5) (a1) {};
		\node[vertex] at (1.5, 2.5) (a2) {};
		\node[special vertex] at (2,   2.5) (a3) {};
		
		\node[special vertex] at (0.5, 2.0) (b1) {};
		\node[vertex] at (0.5, 1.5) (b2) {};
		\node[special vertex] at (0.5, 1.0) (b3) {};
		
		\node[vertex] at (2.5, 2.0) (c1) {};
		\node[special vertex] at (2.5, 1.5) (c2) {};
		\node[vertex] at (2.5, 1.0) (c3) {};
		
		\node[vertex] at (1,   0.5) (d1) {};
		\node[special vertex] at (1.5, 0.5) (d2) {};
		\node[vertex] at (2,   0.5) (d3) {};
		\node (r1) at (1.5,3.0) [] {\faAndroid$_{_1}$};
		\node (r2) at (0.15,  1.5) [] {\faAndroid$_{_2}$};
		\node (r3) at (2.90, 1.5) [] {\faAndroid$_{_3}$};
		\node (r4) at (1.5, 0) [] {\faAndroid$_{_4}$};
		\draw[](a1) -- (a2);
		\draw[](a2) -- (a3);
		\draw[](b1) -- (b2);
		\draw[](b2) -- (b3);
		\draw[](c1) -- (c2);
		\draw[](c2) -- (c3);
		\draw[](d1) -- (d2);
		\draw[](d2) -- (d3);
		
		\draw[](a2) -- (b2);
		\draw[](a2) -- (c1);
		\draw[](c1) -- (d3);
		\draw[](a2) -- (c3);
		\draw[](b2) -- (d1);
		\draw[](b2) -- (d3);
		\begin{pgfonlayer}{background}
			\node[fit=(a1)(a2)(a3),rounded corners,fill=violet!15,inner xsep=3pt,inner ysep=4pt] {};
			\node[fit=(b1)(b2)(b3),rounded corners,fill=green!15,inner xsep=3pt,inner ysep=4pt] {};
			\node[fit=(c1)(c2)(c3),rounded corners,fill=cyan!15,inner xsep=3pt,inner ysep=4pt] {};
			\node[fit=(d1)(d2)(d3),rounded corners,fill=orange!15,inner xsep=3pt,inner ysep=4pt] {};
		\end{pgfonlayer}
		\end{tikzpicture}
		\captionsetup{justification=centering}
		\caption{\small Pose graph $G$}
		\label{fig:pose_graph}
	\end{subfigure}
	~
	\begin{subfigure}[t]{0.20\textwidth}
		\centering
		\begin{tikzpicture}[scale=0.8]
		\tikzstyle{vertex}=[circle,scale=1,draw]
		\node[vertex,fill=violet!15] at (1.5,3.0) (r1) {$x_1$};
		\node[vertex,fill=green!15]  at (0.15,  1.5) (r2) {$x_2$};
		\node[vertex,fill=cyan!15] at (2.90, 1.5) (r3) {$x_3$};
		\node[vertex,fill=orange!15]  at (1.5, 0)  (r4) {$x_4$};
		\draw[thick] (r1) -- (r2);
		\draw[thick] (r1) -- (r3);
		\draw[thick] (r2) -- (r4);
		\draw[thick] (r3) -- (r4);
		\end{tikzpicture}
		\captionsetup{justification=centering}
		\caption{\small Robot-level graph $G_{\Rcal}$}
		\label{fig:robot_graph}
	\end{subfigure}
	\caption
	{\small 
	(a) Example pose graph $G$ with $4$ robots, each with $3$ poses. 
	Each edge denotes a relative pose measurement.
	Private poses are colored in gray. 
	(b) Corresponding robot-level graph $G_{\Rcal}$. 
	Two robots are connected if they share any relative measurements (inter-robot loop closures). 
	Note that at any time during distributed optimization, 
	robots do not need to share their private poses with any other robots.
  	}
	\label{fig:example_graph}
\end{figure}
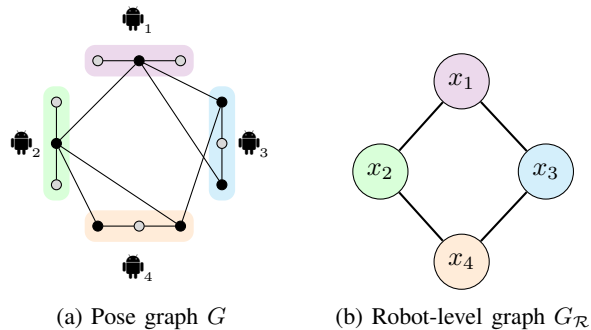

In this section, we formally define pose graph optimization (PGO) in the context of multi-robot SLAM. 
Given relative pose measurements (possibly between different robots), we aim to \emph{jointly} estimate the trajectories of all robots in a global reference frame.
Let $\Rcal = \{1, 2, \hdots, n\}$ be the set of indices associated with $n$ robots. 
Denote the pose of robot $i \in \Rcal$ at time step $\tau$ as $T_{i_\tau} = (R_{i_\tau}, t_{i_\tau}) \in \SE(d)$, where $d \in \{2,3\}$ is the dimension of the estimation problem. 
Here $R_{i_\tau} \in\text{SO}(d)$ is a rotation matrix, and $t_{i_\tau} \in\mathbb{R}^d$ is a translation vector.
A relative pose measurement from $T_{i_\tau}$ to $T_{j_s}$ is denoted as 
$\widetilde{T}^{i_\tau}_{j_s} = (\widetilde{R}^{i_\tau}_{j_s}, \widetilde{t}^{i_\tau}_{j_s}) \in \SE(d)$. 
\blue{We assume the following standard noise model,
which is used by SE-Sync \cite{RosenSESyncIJRR19} and subsequent works 
\cite{Briales2017RAL,Tian2019RBCD},}
\begin{align}
    \widetilde{R}^{i_\tau}_{j_s} &= \underline{R}^{i_\tau}_{j_s} R^{\epsilon}_{i_\tau, j_s},   && R^{\epsilon}_{i_\tau, j_s} \sim \Langevin(I_d, w_R),
    \label{eq:rotation_noise_model} \\
    \widetilde{t}^{i_\tau}_{j_s} &= \underline{t}^{i_\tau}_{j_s} + t^{\epsilon}_{i_\tau, j_s}, && t^{\epsilon}_{i_\tau, j_s} \sim \Ncal(0, w_t^{-1} I_d).
    \label{eq:translation_noise_model}
\end{align}
Above, $\underline{T}^{i_\tau}_{j_s} \!\! = \!\! (\underline{R}^{i_\tau}_{j_s}, \underline{t}^{i_\tau}_{j_s}) \in \SE(d)$ denotes the true (i.e., noiseless) relative transformation.  
The isotropic Langevin noise on rotations \cite{RosenSESyncIJRR19} plays an analogous role as the Gaussian noise on translations. 
\remove{
We note that our formulation trivially generalizes to the case where the values of $w_R$ and $w_t$ vary for different measurements.
In the following, we drop this variation for notation simplicity. }
Given noisy measurements of the form \eqref{eq:rotation_noise_model}-\eqref{eq:translation_noise_model},
we seek to find the maximum likelihood estimates (MLE) for the trajectories of all robots in  $\Rcal$.
Doing so amounts to the following {non-convex} program \cite{RosenSESyncIJRR19}.
\begin{problem}[Maximum Likelihood Estimation]
	\normalfont
	\begin{subequations}
		\begin{align}
		& \underset{}{\text{min}}\!\!\!\!\!
		\sum_{(i_\tau,j_s) \in E} \!\!\!\!\!w_R \norm{R_{j_s} - R_{i_{\tau}} \widetilde{R}^{i_{\tau}}_{j_s}}_\text{F}^2 
		+ w_t \norm{t_{j_s} - t_{i_{\tau}} - R_{i_{\tau}} \widetilde{t}^{i_{\tau}}_{j_s}}^2_2, \nonumber \\
		& \text{s.t.} \;\;\;
		R_{i_{\tau}} \in \SOd(d), \; t_{i_{\tau}} \in \mathbb{R}^d, \; \forall i \in \Rcal, \; \forall \tau.
		\label{eq:mle}
		\tag{$\text{P}_1$}
		\end{align}
	\end{subequations}
	\label{prob:mle}
\end{problem}
\vspace{-0.5cm}

Problem~\eqref{eq:mle} can be compactly represented with a \emph{pose graph} $G \triangleq (V, E)$, where each vertex in $V$ corresponds to a single pose owned by a robot. Observe that the sum in the objective is taken over all edges in $E$, where an edge from $i_{\tau}$ to $j_s$ is formed if there is a relative measurement from $T_{i_{\tau}}$ to $T_{j_s}$. 
Fig.~\ref{fig:pose_graph} shows an example pose graph.

In this paper, we further consider the \emph{rank-restricted relaxation} of \eqref{eq:mle}, \blue{first proposed by SE-Sync \cite{RosenSESyncIJRR19}.
}
Denote the Stiefel manifold as $\Stiefel(d,r) \triangleq \{Y \in \mathbb{R}^{r \times d}: Y^\top Y = I_d\}$, where $r \geq d$.
The rank-$r$ relaxation of \eqref{eq:mle} is defined as the following non-convex Riemannian optimization problem.
\begin{problem}[Rank-$r$ Relaxation]
	\normalfont
	\begin{subequations}
		\begin{align}
		& \underset{}{\text{min}}\!\!\!\!\!
		\sum_{(i_{\tau},j_s) \in E} \!\!\!\! w_R \norm{Y_{j_s} - Y_{i_{\tau}} \widetilde{R}^{i_{\tau}}_{j_s}}_\text{F}^2 
		+ w_t \norm{p_{j_s} - p_{i_{\tau}} - Y_{i_{\tau}} \widetilde{t}^{i_{\tau}}_{j_s}}^2_2, \nonumber \\
		& \text{s.t.} \;\;\;
		Y_{i_{\tau}} \in \Stiefel(d,r), \; p_{i_{\tau}} \in \mathbb{R}^r, \; \forall i \in \Rcal, \; \forall \tau.
		\tag{$\text{P}_2$}
		\label{eq:qcqp}
		\end{align}
	\end{subequations}
	\label{prob:qcqp}
\end{problem}
\vspace{-0.5cm}

Observe that for $r = d$, the Stiefel manifold is identical to the orthogonal group $\Stiefel(d,d) = \Orthogonal(d)$. In this case, \eqref{eq:qcqp} is referred to as the \emph{orthogonal relaxation} of \eqref{eq:mle}, obtained by dropping the determinant constraint on $\SOd(d)$. 
As $r$ increases beyond $d$, we obtain a hierarchy of rank-restricted problems, each having the form of \eqref{eq:qcqp} but with a slightly ``lifted'' search space as determined by $r$.
\blue{
The idea of rank relaxations is first proposed in SE-Sync \cite{RosenSESyncIJRR19}, where Rosen et al. use a very similar hierarchy of rank relaxations as a proxy to solve the semidefinite relaxation of PGO, thereby recovering \emph{global minimizers} to the original MLE \eqref{eq:mle}.\footnote{
\blue{
The only difference between the rank relaxations in \cite{RosenSESyncIJRR19} and \eqref{eq:qcqp} is that translations variables are analytically eliminated in \cite{RosenSESyncIJRR19}.
In \eqref{eq:qcqp}, we keep translations in our optimization to retain sparsity. 
This approach is first used in Cartan-Sync~\cite{Briales2017RAL}.}
}
This elegant result motivates us to consider \eqref{eq:qcqp} in addition to the original MLE formulation \eqref{eq:mle}.
While SE-Sync assumes a centralized (single-robot) scenario, we focus on distributed PGO and develop an \emph{asynchronous} algorithm to solve both \eqref{eq:mle} and \eqref{eq:qcqp} in the presence of communication delay.}

\remove{Once we solve \eqref{eq:qcqp}, either globally or locally to a critical point, we can apply a distributed \emph{rounding} procedure (e.g., as detailed in \cite{Tian2019RBCD}) to obtain a feasible solution to the original MLE problem \eqref{eq:mle}.
In addition, note that \eqref{eq:qcqp} shares the same sparsity structure as encoded by the pose graph.}

For the purpose of designing decentralized algorithms (Section~\ref{sec:algorithm}), it is more convenient to rewrite \eqref{eq:mle} and \eqref{eq:qcqp} into a more abstract form at the level of robots, which may be done as follows.
\begin{problem}[Robot-level Optimization Problem]
	\normalfont
	\begin{equation}
		\begin{aligned}
		& \underset{}{\text{min}}
		\sum_{(i,j) \in E_{\Rcal}} f_{ij}(x_i, x_j) + \sum_{i \in \Rcal} h_i(x_i),  \\
		& \text{s.t.} \;\;\;
		x_i \in \Mcal_i, \; \forall i \in \Rcal.
		\end{aligned}
		\tag{$\text{P}$}
		\label{eq:P}
	\end{equation}
	\label{prob:P}
\end{problem}
\vspace{-0.5cm}

In \eqref{eq:P}, each variable $x_i$ concatenates all variables owned by robot $i \in \Rcal$. 
For instance, for \eqref{eq:qcqp}, $x_i$ contains all the ``lifted'' rotation  and translation variables of robot $i$. 
Let $n_i$ be the number of poses of robot $i$.
Then,
\begin{align}
    x_i &= \begin{bmatrix}
    Y_{i_1} & p_{i_1} & \hdots & Y_{i_{n_i}} & p_{i_{n_i}}
    \end{bmatrix}, \\
    \Mcal_i &= (\Stiefel(d,r) \times \mathbb{R}^r)^{n_i}.
\end{align}

The cost function in \eqref{eq:P} consists of a set of \emph{shared costs} $f_{ij}: \Mcal_i \times \Mcal_j \to \mathbb{R}$ between pairs of robots, and a set of \emph{private costs} $h_i: \Mcal_i \to \mathbb{R}$ for individual robots.  
Intuitively, $f_{ij}$ is formed by relative measurements between any of robot $i$'s poses and $j$'s poses. 
In contrast, $h_i$ is formed by relative measurements within robot $i$'s own trajectory. 

Similar to the way a pose graph is defined, we can encode the structure of \eqref{eq:P} using a \emph{robot-level graph} $G_{\Rcal} \triangleq (\Rcal, E_{\Rcal})$; see Fig.~\ref{fig:robot_graph}.
$G_{\Rcal}$ can be viewed as a ``reduced'' graph of the pose graph, in which each vertex corresponds to the entire trajectory of a single robot $i \in \Rcal$.
Two robots $i,j$ are connected in $G_{\Rcal}$ if they share any relative measurements $(i_{\tau}, j_s) \in E$.
In this case, we call $j$ a \emph{neighboring robot} of $i$, 
and $j_s$ a \emph{neighboring pose} of robot $i$.
If a pose variable is not a neighboring pose to any other robots, we call this pose a \emph{private pose} \cite{Tian2019RBCD}.
We note that for robot $i$ to evaluate the shared cost $f_{ij}$, it only needs to know its neighboring poses in robot $j$'s trajectory (see Fig.~\ref{fig:example_graph}). 
This property is crucial in preserving the \emph{privacy} of participating robots \cite{Choudhary2017Distributed,Tian2019RBCD}, i.e., at any time, a robot does not need to share its private poses with any of its teammates.

\section{Proposed Algorithm}
\label{sec:algorithm}

We present our main algorithm, Asynchronous Stochastic Parallel Pose Graph Optimization ($\ALGORITHM$), for solving distributed PGO problems of the form \eqref{eq:P}.
Our algorithm is inspired by asynchronous stochastic coordinate descent (e.g., see~\cite{Liu2015SIAM}).
In the context of distributed PGO, each coordinate corresponds to the stacked relative pose observations $x_i$ of a single robot as defined in \eqref{eq:P}. 

In a practical multi-robot SLAM scenario, each robot can optimize its own pose estimates at any time, and can additionally share its (non-private) poses with others when communication is available. 
Correspondingly, each robot running $\ALGORITHM$ has two concurrent onboard processes, which we refer to as the \emph{optimization thread} and \emph{communication thread}.
We emphasize that the robots perform both optimization and communication completely in parallel and without synchronization with each other.
We begin by describing the communication thread and then proceed to the optimization thread. 
Without loss of generality, we describe the algorithm from the perspective of robot $i \in \Rcal$.

\subsection{Communication Thread}
\label{subsec:communication_thread}

As part of the communication module, each robot $i \in \Rcal$ implements a local data structure, called a \emph{cache}, that contains the robot's own variable $x_i$, together with the most recent copies of neighboring poses received from the robot's neighbors. 
\blue{A very similar design that allows asynchronous communication is proposed by Cunningham et al.~\cite{CunninghamDDFSAM,Cunningham2012ICRA,CunninghamDDFSAM2}, although the authors have not discussed convergence in the asynchronous setting.}

Since only $i$ can modify $x_i$, the value of $x_i$ in robot $i$'s cache is guaranteed to be up-to-date at anytime. 
In contrast, the copies of neighboring poses from other robots can be \emph{out-of-date} due to communication delay. 
For example, by the time robot $i$ receives and uses a copy of robot $j$'s poses, $j$ might have already updated its poses due to its local optimization process. 
In Section~\ref{sec:convergence}, we show that $\ALGORITHM$ is resilient against such network delay.
Nevertheless, for $\ALGORITHM$ to converge, we still assume that the total delay induced by the communication process remains \emph{bounded} (Section~\ref{sec:convergence}).
The communication thread performs the following two operations over the cache.

\vspace{0.1cm}
\noindent $\bullet$ \textbf{Receive}:
After receiving a neighboring pose, e.g., $(R_{j_s}, t_{j_s})$ from a neighboring robot $j$ over the network, the communication thread updates the corresponding entry in the cache to store the new value. 

\vspace{0.1cm}
\noindent $\bullet$ \textbf{Send}:
Periodically (when communication is available), robot $i$ also transmits its latest public pose variables (i.e., poses that have inter-robot measurements with other robots) to its neighbors.
Recall from Section~\ref{sec:problem_formulation} that robot $i$ does not need to send its private poses, as these poses are not needed by other robots to optimize their estimates.

\subsection{Optimization Thread}
\label{subsec:optimization_thread}
Concurrent to the communication thread, the optimization thread is invoked by a local clock that ticks according to a Poisson process of rate $\lambda > 0$.

\begin{definition}[Poisson process \cite{Tijms2004book}]
Consider a sequence $\{X_1, X_2,...\}$ of positive, independent random variables that represent the time elapsed between consecutive events (in this case, clock ticks). 
Let $N(t)$ be the number of events up to time $t \geq 0$.
The counting process $\{N(t), t \geq 0\}$ is a Poisson process with rate $\lambda > 0$ if the interarrival times $\{X_1, X_2,...\}$ have a common exponential distribution function, 
\begin{equation}
    P(X_k \leq a) = 1 - e^{-\lambda a}, \; a \geq 0.
\end{equation}
\end{definition}

The use of Poisson clocks originates from the design of randomized gossip algorithms by Boyd et al.~\cite{Boyd2006Gossip} and is a commonly used tool for analyzing the global behavior of distributed randomized algorithms.
The rate parameter $\lambda$ is equal among robots.
In practice, we can adjust $\lambda$ based on the extent of network delay and the robots' computational capacity.
Using this local clock, the optimization thread performs the following operations in a loop.

\vspace{0.2cm}
\noindent $\bullet$ \textbf{Read}:
For each neighboring robot $j \in \Nbr(i)$, read the value of $x_j$ stored in the local cache. Denote the read values as $\hat{x}_j$. 
Recall that $\hat{x}_j$ can be \emph{outdated}, for example if robot $i$ has not received the latest messages from robot $j$.
In addition, read the value of $x_i$, denoted as $\hat{x}_i$.
Recall from Section~\ref{subsec:communication_thread} that $\hat{x}_i$ is guaranteed to be up-to-date.

In practice, $\hat{x}_j$ only contains the set of neighboring poses from robot $j$ since $f_{ij}$ is independent from the rest of $j$'s poses (Fig.~\ref{fig:example_graph}). 
However, for ease of notation and analysis, we treat $\hat{x}_j$ as if it contains the entire set of $j$'s poses.

\vspace{0.2cm}
\noindent $\bullet$ \textbf{Compute}:
Form the local cost function for robot $i$, denoted as $g_i(x_i): \Mcal_i \to \mathbb{R}$, by aggregating relevant costs in \eqref{eq:P} that involve $x_i$, 
\begin{equation}
    g_i(x_i) =  h_i(x_i)  + \sum_{j \in \Nbr(i)} f_{ij}(x_i, \hat{x}_j).
    \label{eq:local_cost}
\end{equation}
Compute the Riemannian gradient at robot $i$'s current estimate $\hat{x}_i$,
\begin{equation}
    \eta_i = \rgrad g_i(\hat{x}_i) \in T_{\hat{x}_i} \Mcal_i.
    \label{eq:local_descent}
\end{equation}

\vspace{0.2cm}
\noindent $\bullet$ \textbf{Update}:
At the next local clock tick, update $x_i$ in the direction of the negative gradient, 
\begin{equation}
    x_i \leftarrow \Retr_{\hat{x}_i}(-\gamma \eta_i),
    \label{eq:update_step}
\end{equation}
where $\gamma > 0$ is a constant stepsize.
Equation~\eqref{eq:update_step} gives the simplest update rule that robots can follow.
\blue{In Section~\ref{sec:convergence}, we further prove convergence for this update rule.

To accelerate convergence in practice, SE-Sync~\cite{RosenSESyncIJRR19} and Cartan-Sync~\cite{Briales2017RAL} use a heuristic known as \emph{preconditioning}, 
which is also applicable to $\ALGORITHM$. 
With preconditioning, the following alternative update direction is used,
}
\begin{equation}
    x_i \leftarrow \Retr_{\hat{x}_i}(-\gamma \Precon g_i(\hat{x}_i)[\eta_i]).
    \label{eq:preconditioned_update_step}
\end{equation}
In~\eqref{eq:preconditioned_update_step}, $\Precon g_i(\hat{x}_i): T_{\hat{x}_i} \Mcal_i \to T_{\hat{x}_i} \Mcal_i$ is a linear, symmetric, and positive definite mapping on the tangent space that approximates the inverse of Riemannian Hessian.
Intuitively, preconditioning helps first-order methods benefit from using the (approximate) second-order geometry of the cost function, which often results in significant speedup especially on poorly conditioned problems.

\subsection{Implementation Details}

To make the local clock model valid, we require that the total execution time of the \textbf{Read}-\textbf{Compute}-\textbf{Update} sequence be smaller than the interarrival time of the Poisson clock, so that the current sequence can finish before the next one starts. 
This requirement is fairly lax in practice, as all three steps only involve minimal computation and access to local memory.
In the worst case, since the interarrival time is determined by $1/\lambda$ on average \cite{Tijms2004book}, one can also decrease the clock rate $\lambda$ to create more time for each update.

In addition, we note that although the optimization and communication threads run concurrently, minimal thread-level synchronization is required to ensure the so-called \emph{atomic read and write} of individual poses.
Specifically, a thread cannot read a pose in the cache if the other thread is actively modifying its value (otherwise the read value would not be valid).
Such synchronization can be easily enforced using software locks. 
\remove{
In practice, however, due to the large number of poses owned by each robot, the aforementioned synchronization only happens relatively rarely. }

\section{Convergence Analysis}
\label{sec:convergence}

\subsection{Global View of the Algorithm}
\label{subsec:global_view}

In Section~\ref{sec:algorithm}, we described $\ALGORITHM$ from the local perspective of each robot. 
For the purpose of establishing convergence, however, we need to analyze the systematic behavior of this algorithm from a global perspective \cite{Liu2015SIAM,Liu2015JMLR,Lian2018ICML,Boyd2006Gossip}.
To do so, let $k = 0, 1, \hdots $ be a virtual counter that counts the total number of \textbf{Update} operations applied by all robots. 
In addition, let the random variable $\ik \in \Rcal$ represent the robot that updates at global iteration $k$. 
We emphasize that $k$ and $\ik$ are purely used for theoretical analysis, 
and are unknown to any of the robots in practice. 

Recall from Section~\ref{subsec:optimization_thread} that all \textbf{Update} steps are generated by $ n = |\Rcal|$ independent Poisson processes, each with rate $\lambda$.  
In the global perspective, merging these local processes is equivalent to creating a single, global Poisson clock with rate $\lambda n$.
Furthermore, at any time, all robots have equal probabilities of generating the next \textbf{Update} step, i.e., for all $k \in \mathbb{N}$, $\ik$ is i.i.d. uniformly distributed over the set $\Rcal$. See \cite{Tijms2004book} for proofs of these results.

\begin{algorithm}[t]
	\caption{{\small \textsc{Global View of $\ALGORITHM$} (For Analysis Only)} }
	\label{alg:global}
	\begin{algorithmic}[1]
		\renewcommand{\algorithmicrequire}{\textbf{Input:}}
		\renewcommand{\algorithmicensure}{\textbf{Output:}}
		\Require
		\Statex Initial solution $x^0 \in \Mcal$ and stepsize $\gamma > 0$.
		\vspace{0.1cm}
		\For{\text{global iteration} $k=0,1,\hdots$}
		\State Select robot $\ik \in \Rcal$ uniformly at random. \label{line:sample}
		\State Read $\hat{x}_{\ik} = x^k_{\ik}$. \label{line:read_a}
		\State Read $\hat{x}_{\jk} = x^{k - B(\jk)}_{\jk}, \; \forall \jk \in \Nbr(\ik)$. \label{line:read_b}
		\State Compute local gradient $\eta^k_{\ik} = \rgrad g_{\ik} (\hat{x}_{\ik})$. \label{line:compute}
		\State Update $x_{\ik}^{k+1} = \Retr_{\hat{x}_{\ik}}(-\gamma \eta^k_{\ik})$. \label{line:update_a}
		\State Carry over all $x_{j}^{k+1} = x_{j}^k, \; \forall j \neq \ik$.
		\label{line:update_b}
		\EndFor
	\end{algorithmic}
\end{algorithm}

Using this result, we can write the iterations of $\ALGORITHM$ from the global view; see Algorithm~\ref{alg:global}.
We use $x^k \triangleq \begin{bmatrix} x^k_1 & x^k_2 & \hdots & x^k_n \end{bmatrix}$ to represent the value of all robots' poses after $k$ global iterations (i.e., after $k$ total \textbf{Update} steps). 
Note that $x$ lives on the product manifold $\Mcal \triangleq \Mcal_1 \times \Mcal_2 \times  \hdots \Mcal_n$. 
At global iteration $k$, a robot $i_k$ is selected from $\Rcal$ uniformly at random (line~\ref{line:sample}). 
Robot $i_k$ then follows the steps in Section~\ref{subsec:optimization_thread} to update its own variable (line~\ref{line:read_a}-\ref{line:update_a}).
We have used the fact that $\hat{x}_{\ik}$ is always up-to-date (line~\ref{line:read_a}), while $\hat{x}_{\jk}$ is outdated for $B(\jk)$ total \textbf{Update} steps (line~\ref{line:read_b}).
Except robot $\ik$, all other robots do not update (line~\ref{line:update_b}).
As an additional notation that will be useful for later analysis, 
we note that line~\ref{line:update_b} can be equivalently written as $x_{j}^{k+1} = \Retr_{x^k_{j}}(-\gamma \eta^k_{j})$ with $\eta^k_{j} = 0$.

\subsection{Sufficient Conditions for Convergence}

We establish sufficient conditions for $\ALGORITHM$ to converge to first-order critical points.
Due to space limitation, all proofs are deferred to the appendix~\cite{Appendix}.
We adopt the commonly used \emph{partially asynchronous} model \cite{bertsekas1989parallel}, which assumes that delay caused by asynchrony is not arbitrarily large. 
In practice, the magnitude of delay is affected by various factors such as the rate of communication (Section~\ref{subsec:communication_thread}), the rate of local optimization (Section~\ref{subsec:optimization_thread}), and intrinsic network latency.
For the purpose of analysis, we assume that all these factors can be summarized into a single constant $B$, which bounds the maximum delay in terms of number of \emph{global iterations} (i.e., \textbf{Update} steps applied by all robots) in Algorithm~\ref{alg:global}.

\begin{assumption}[Bounded Delay]
In Algorithm~\ref{alg:global}, there exists a constant $B \geq 0$ such that $B(\jk) \leq B$ for all $\jk$.
\label{as:bounded_delay}
\end{assumption}
\blue{
Assumption~\ref{as:bounded_delay} imposes a worst-case upper bound on the delay, and allows the actual delay to fluctuate within this upper bound.}
In addition, for both the MLE problem \eqref{eq:mle} and its rank-restricted relaxations \eqref{eq:qcqp}, the gradients enjoy a Lipschitz-type condition, which is proved in our previous work \cite{Tian2019RBCD} and will be used extensively in the rest of the analysis. 

\begin{lemma}[Lipschitz-type gradient for pullbacks \cite{Boumal2018Convergence,Tian2019RBCD}]
Denote the cost function of \eqref{eq:mle} and \eqref{eq:qcqp} as $f: \Mcal \to \mathbb{R}$.
Define the \emph{pullback} cost as $\hat{f}_x \triangleq f \circ \Retr_x: T_x \Mcal \to \mathbb{R}$. 
There exists a constant $L \geq 0$ such that for any $x \in \Mcal$ and $\eta \in T_x \Mcal$,
\begin{equation}
\big|\fhat_x(\eta) - [f(x) + \langle \eta, \rgrad_x f \rangle]\big| \leq \frac{L}{2} \norm{\eta}^2.
\label{eq:lipschitz_gradient_pullback}
\end{equation} 
\label{lem:lipschitz_gradient_pullback}
\end{lemma}

\vspace{-0.5cm}
The condition \eqref{eq:lipschitz_gradient_pullback} is first proposed by \cite{Boumal2018Convergence} as an adaptation of Lipschitz continuous gradient to Riemannian optimization. 
Using the bounded delay assumption and the Lipschitz-type condition in \eqref{eq:lipschitz_gradient_pullback}, we can proceed to analyze the change in cost function after a single iteration of Algorithm~\ref{alg:global} (in the global view). 
We formally state the result in the following lemma. 

\begin{lemma}[Descent Property of Algorithm~\ref{alg:global}]
Under Assumption~\ref{as:bounded_delay}, each iteration of Algorithm~\ref{alg:global} satisfies,
\begin{equation}
\begin{aligned}
    f(x^{k+1}) - f(x^k) \leq
    -\frac{\gamma}{2} \norm{\rgrad_{\ik} f(x^k)}^2
    - \frac{\gamma-L\gamma^2}{2} \norm{\eta^k_{\ik}}^2 \\
    + \frac{\Delta B L^2 \alpha^2 \gamma^3 }{2} 
    \sum_{\jk \in \Nbr(\ik)} \sum_{k' = k - B}^{k-1} \norm{\eta_{\jk}^{k'}}^2,
    \label{eq:descent_lemma}
\end{aligned}
\end{equation}
where $\eta_i^k$ denotes the update taken by robot $i$ at iteration $k$,
$\alpha > 0$ is a constant related to the retraction, and $\Delta > 0$ is the maximum degree of the robot-level graph $G_{\Rcal}$.
\label{lem:descent_lemma}
\end{lemma}
In \eqref{eq:descent_lemma}, the last term on the right hand side sums over the squared norms of a set of $\{\eta_{\jk}^{k'}\}$, where each $\eta_{\jk}^{k'}$ corresponds to the update taken by a neighbor $\jk$ at an earlier iteration $k'$.
This term is a direct consequence of delay in the system, 
and is also the main obstacle for proving convergence in the asynchronous setting. 
Indeed, without this term, it is straightforward to verify that any stepsize that satisfies $0 < \gamma < 1/L$ guarantees $f(x^{k+1}) \leq f(x^k)$, and thus leads to convergent behavior. 
With the last term in \eqref{eq:descent_lemma}, however, the overall cost could increase after each iteration.

While the delay-dependent error term gives rise to additional challenges, our next theorem states that with sufficiently small stepsize, this error term is {inconsequential} and $\ALGORITHM$ provably converges to first-order critical points.

\begin{theorem}[Global convergence of $\ALGORITHM$]
Let $f^\star$ be any global lower bound on the optimum of \eqref{eq:P}. 
Define $\rho \triangleq \Delta / n$. 
Let $\bar{\gamma} > 0$ be an upper bound on the stepsize that satisfies,
\begin{equation}
    2 \rho \alpha^2 B^2 L^2 \bar{\gamma}^2 + L \bar{\gamma} - 1 \leq 0.
    \label{eq:stepsize_inequality}
\end{equation}
In particular, the following choice of $\bar{\gamma}$ satisfies \eqref{eq:stepsize_inequality}:
\begin{equation}
    \bar{\gamma} =
    \begin{cases}
    \frac
    {\sqrt{1 + 8 \rho \alpha^2 B^2 } - 1}
    {4 \rho \alpha^2 B^2 L}, & B > 0, \\
    1 / L, & B = 0.
    \end{cases}
    \label{eq:stepsize_bound}
\end{equation}
Under Assumption~\ref{as:bounded_delay},
if $0 < \gamma \leq \bar{\gamma}$, $\ALGORITHM$ converges to a first-order critical 
point with global sublinear rate. Specifically, after $K$ total update steps, 
\begin{equation}
    \min_{k \in \{0,\hdots,K-1\}} \mathbb{E}_{i_{0:K-1}} \bigg[\norm{\rgrad f(x^k)}^2 \bigg] \leq 
    \frac{2n(f(x^0) - f^\star)}{\gamma K}.
    \label{eq:exact_convergence}
\end{equation}
\label{thm:exact_convergence}
\end{theorem}
\vspace{-0.5cm}

\begin{remark}
\normalfont

To the best of our knowledge, Theorem~\ref{thm:exact_convergence} establishes the first convergence result for asynchronous algorithms when solving a \emph{non-convex} optimization problem over the product of matrix manifolds.  
While the existence of a convergent stepsize $\bar{\gamma}$ is of theoretical importance, 
we further note that its expression $\eqref{eq:stepsize_bound}$ offers the correct qualitative insights with respect to various problem-specific parameters, which we discuss next. 

\vspace{0.2cm}
\noindent \underline{\textit{Relation with maximum delay ($B$)}}: 
$\bar{\gamma}$ increases as maximum delay $B$ decreases. 
Intuitively, as communication becomes increasingly available, each robot may take larger steps without causing divergence. 
The inverse relationship between $\bar{\gamma}$ and $B$ is well known in the asynchronous optimization literature, and is first established by Bertsekas and Tsitsilis~\cite{bertsekas1989parallel} in the Euclidean setting.

\vspace{0.2cm}
\noindent \underline{\textit{Relation with problem sparsity ($\rho$)}}: 
$\bar{\gamma}$ increases as $\rho$ decreases.
Recall that $\rho \triangleq \Delta / n$ is defined as the ratio between the maximum number of neighbors a robot has and the total number of robots.
Thus, $\rho$ is a measure of \emph{sparsity} of the robot-level graph $G_{\Rcal}$. 
Intuitively, as $G_{\Rcal}$ becomes more sparse, robots can use larger stepsize as their problems become increasingly decoupled.
Such positive correlation between $\bar{\gamma}$ and problem sparsity has been a crucial feature in state-of-the-art asynchronous algorithms; see e.g., \cite{Niu2011NIPS}. 

\vspace{0.2cm}
\noindent \underline{\textit{Relation with problem smoothness ($L$)}}: 
From \eqref{eq:stepsize_bound}, it can been seen that $\bar{\gamma}$ increases asymptotically with $\mathcal{O}(1/L)$.
Moreover, when there is no delay ($B = 0$), our stepsize matches the well-known constant of $1/L$ with which synchronous gradient descent converges to first-order critical points; see e.g., \cite{Boumal2018Convergence}. 

\end{remark}


\section{Experimental Results}
\label{sec:experiments}

We implement $\ALGORITHM$ in C++ and evaluate its performance on both simulated and real-world PGO datasets.
We use ROPTLIB~\cite{ROPTLIB} for manifold related computations, and the Robot Operating System (ROS)~\cite{ROS} for inter-robot communication.  
The Poisson clock is implemented by halting the optimization thread after each iteration for a random amount of time exponentially distributed with rate $\lambda$ (default to $1000$~Hz). 
Since the time taken by each iteration is negligible, we expect the practical difference between this implementation and the theoretical model in Section~\ref{subsec:optimization_thread} to be insignificant.
All robots are simulated as separate ROS nodes running on a desktop computer with an Intel i7 quad-core CPU and $16$~GB memory. 


For each PGO problem, we use $\ALGORITHM$ to solve its rank-restricted relaxation $\eqref{eq:qcqp}$ with $r = 5$.
As is commonly done in prior work \cite{Agarwal2011NIPS,Niu2011NIPS,Liu2015SIAM,Liu2015JMLR,Lian2015NIPS,Lian2018ICML,Cannelli2019},
in our experiments we select the stepsize empirically.
During optimization, we record the evolution of the Riemannian gradient norm $\norm{\rgrad f(x^k)}$, which measures convergence to a first-order critical point.
In addition, we also record the optimality gap $f(x^k) - f(x^\star)$, where $x^\star$ is a global minimizer to the PGO problem \eqref{eq:mle} computed by Cartan-Sync \cite{Briales2017RAL}.
In Section~\ref{subsec:benchmark_results}, we also round the solution to $\SE(d)$ using the method in \cite{Tian2019RBCD} and then compute the 
\blue{translation root mean squared error (RMSE) and rotation RMSE (in chordal distance) with repspect to the global minimizer.
}

\subsection{Evaluation in Simulation}
\label{subsec:simulation_results}

\begin{figure*}[t]
	\centering
	\begin{subfigure}[t]{0.3\textwidth}
		\centering
		\includegraphics[trim= 90 0 50 90, width=\textwidth]{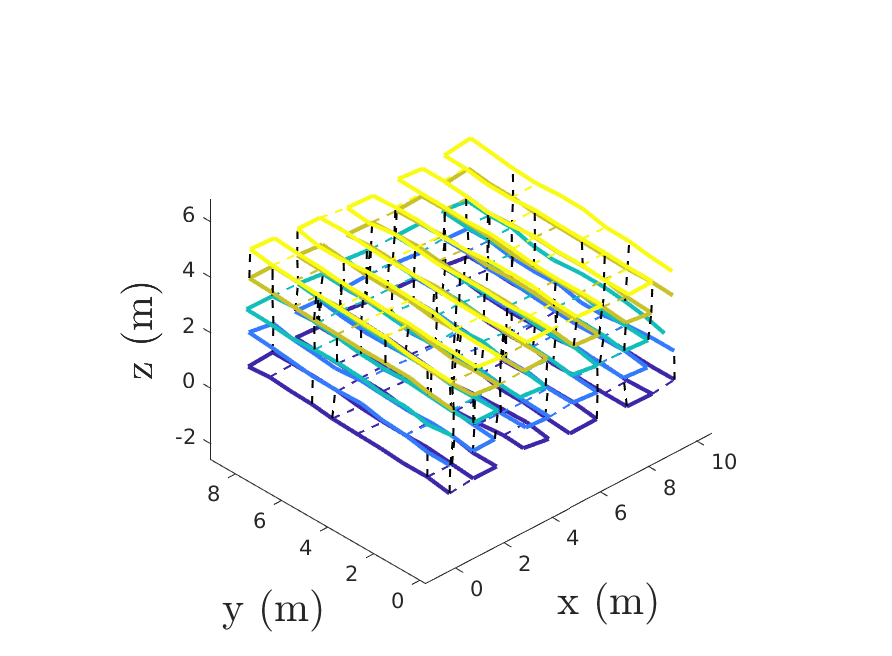}
		\caption{\small Simulation}
		\label{fig:simulation_solution}
	\end{subfigure}
	~
	\begin{subfigure}[t]{0.3\textwidth}
		\centering
		\includegraphics[width=\textwidth]{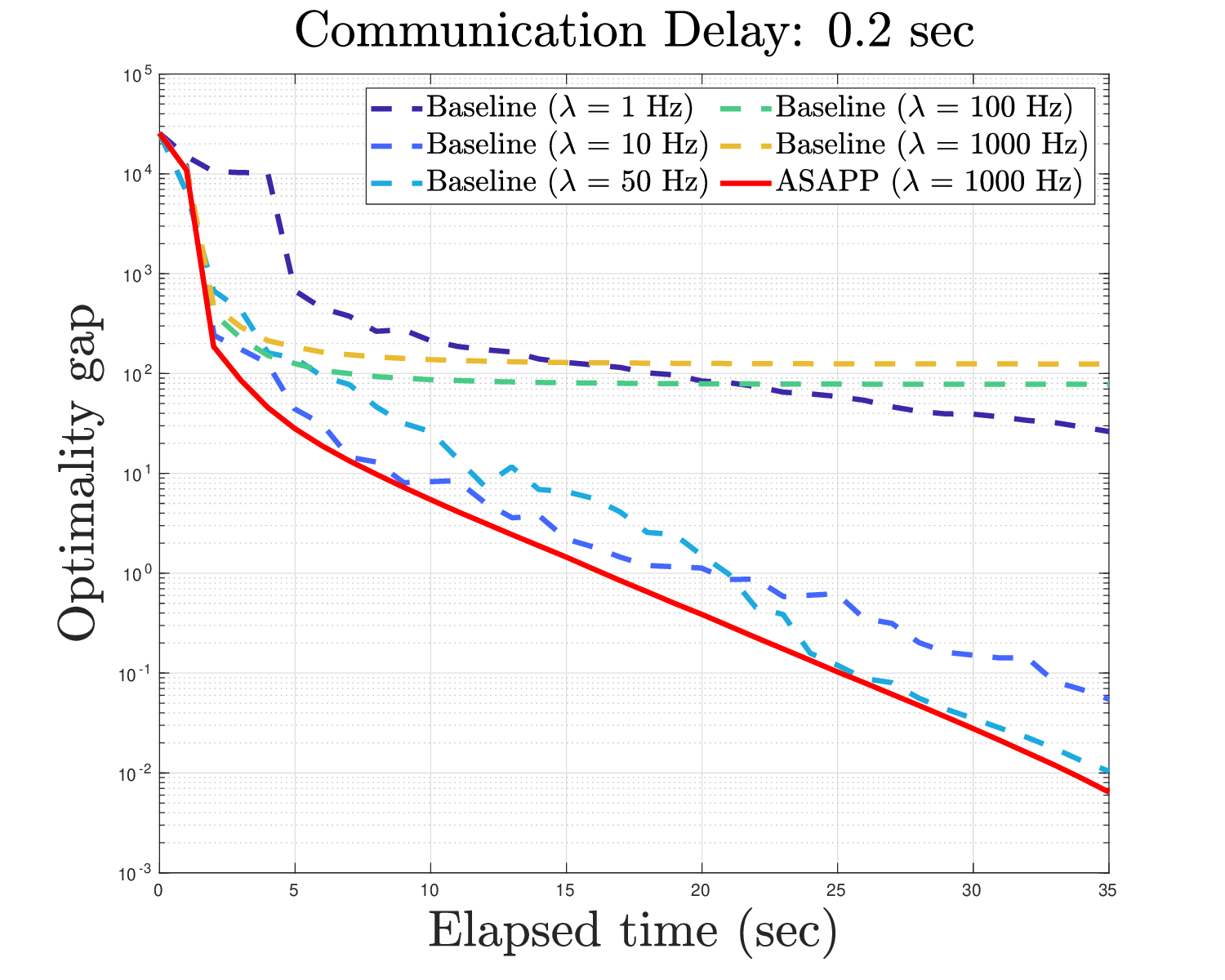}
		\caption{\small Optimality gap}
		\label{fig:simulation_optimality}
	\end{subfigure}
	~
	\begin{subfigure}[t]{0.3\textwidth}
		\centering
		\includegraphics[width=\textwidth]{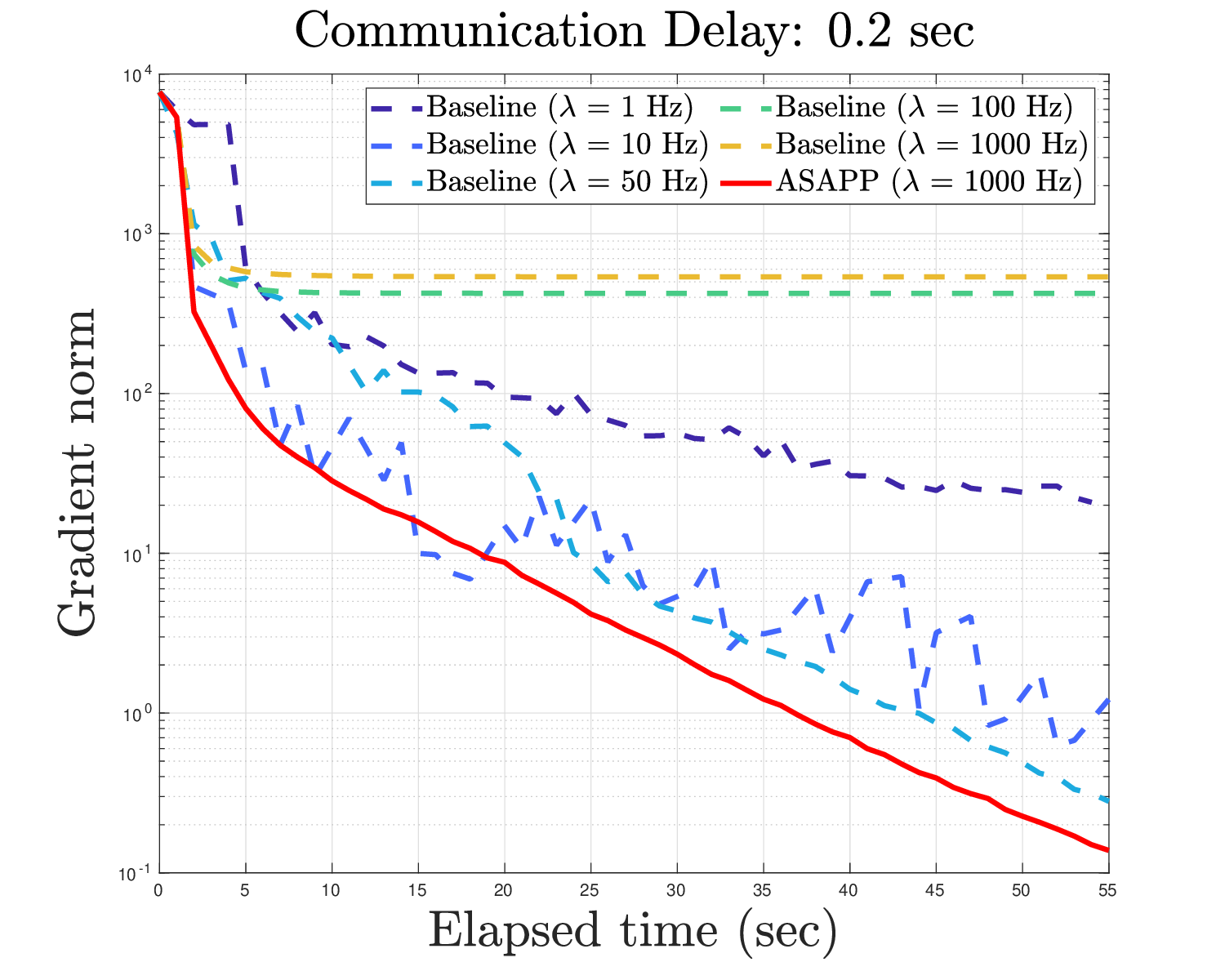}
		\caption{\small Riemannian gradient norm}
		\label{fig:simulation_gradnorm}
	\end{subfigure}
	\caption{\small 
	    Performance evaluation on 5 robot simulation.
	    The communication delay is fixed at $0.5$~s.
	    We compare $\ALGORITHM$ (with stepsize $\gamma = 5\times 10^{-4}$) with a baseline algorithm in which each robot uses Riemannian trust-region method to optimize its local variables.
	    For a comprehensive evaluation, we run the baseline with varying optimization rate to record its performance under both synchronous and asynchronous regimes.
	    (a) Example trajectories estimated by $\ALGORITHM$, where trajectories of 5 robots are shown in different colors. 
	    Inter-robot measurements (loop closures) are shown as black dashed lines.
		(b) Optimality gap with respect to the centralized global minimizer $f(x^k) - f(x^\star)$.
		(c) Riemannian gradient norm $\norm{\rgrad f(x^k)}$.
		}
	\label{fig:simulation_experiment}
\end{figure*} 

\begin{figure}[t]
	\centering
	\includegraphics[width=0.3\textwidth]{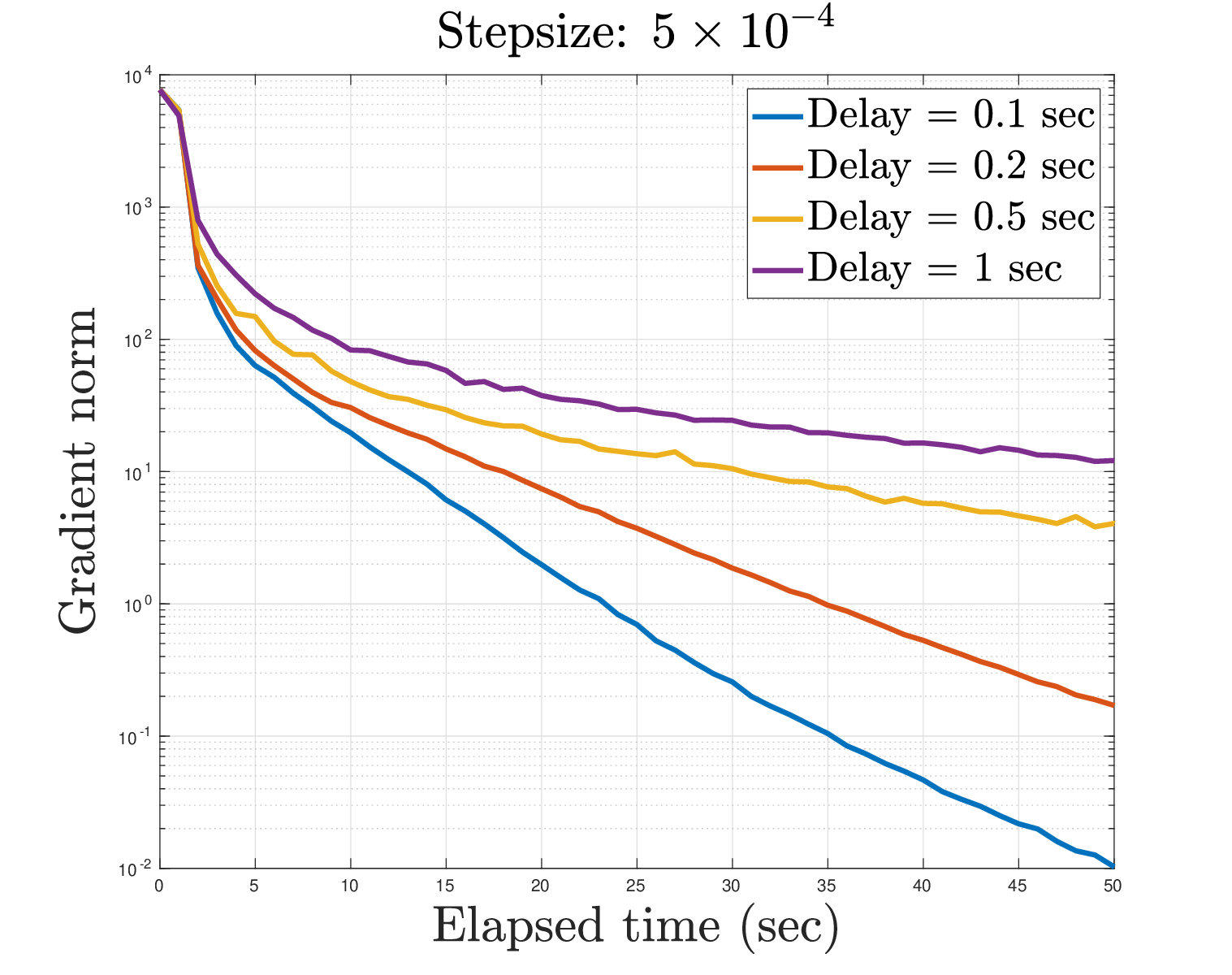}
	\caption{\small 
	    Convergence speed of $\ALGORITHM$ (stepsize $\gamma = 5 \times 10^{-4}$) with varying communication delay. 
	    As delay decreases, convergence becomes faster because robots have access to more up-to-date information from each other. 
	    }
	\label{fig:simulation_vary_delay}
\end{figure}

We evaluate $\ALGORITHM$ in a simulated multi-robot SLAM scenario in which $5$ robots move next to each other in a 3D grid with lawn mower trajectories (Fig.~\ref{fig:simulation_solution}). 
Each robot has $100$ poses. 
With probability $0.3$, loop closures within and across trajectories are generated for poses within $1$~m of each other.
All measurements are corrupted by Langevin rotation noise with $2^\circ$ standard deviation, and Gaussian translation noise with $0.05$~m standard deviation. 
\blue{To minimize communication during initialization, we initialize the solution by propagating relative measurements along a spanning tree of the global pose graph.}
The stepsize used in simulation is $\gamma = 5\! \times\! 10^{-4}$.

In the first experiment, we simulate communication delay by letting each robot communicate every $0.2$~s.
We compare the performance of $\ALGORITHM$ (without preconditioning) against a baseline algorithm in which each robot uses the second-order Riemannian trust-region (RTR) method to optimize its local variable, similar to the approach in \cite{Tian2019RBCD}. 
\blue{Starting with SE-Sync \cite{RosenSESyncIJRR19}, RTR has been used as the default solver in centralized or synchronous settings due to its global convergence guarantees and ability to exploit second-order geometry of the cost function.}
For a comprehensive evaluation, we record the performance of this baseline at different optimization rates (i.e. frequency at which robots update their local trajectories).

Fig.~\ref{fig:simulation_optimality} shows the optimality gaps achieved by the evaluated algorithms as a function of wall clock time.
The corresponding reduction in the Riemannian gradient norm is shown in Fig.~\ref{fig:simulation_gradnorm}.
$\ALGORITHM$ outperforms all variants of the baseline algorithm (dashed curves). 
We note that the behavior of the baseline algorithm is expected. 
At a low rate, e.g., $\lambda = 1$~Hz (dark blue dashed curve), the baseline algorithm is essentially synchronous as each robot has access to up-to-date poses from others. 
The empirical convergence speed is nevertheless slow, since each robot needs to wait for up-to-date information to arrive after each iteration.
At a high rate, e.g., $\lambda = 1000$~Hz (dark yellow dashed curve), robots essentially behave asynchronously. 
However, since RTR does not regulate stepsize at each iteration, robots often significantly alter their solutions in the wrong direction (as a result of using outdated information), which leads to slow convergence or even non-convergence. 
\remove{
Lastly, we observe that at an intermediate rate, e.g., $\lambda = 50$~Hz, convergence speed of the baseline approaches that of $\ALGORITHM$. 
However, we emphasize that the baseline algorithm does not provide any convergence guarantees.}
In contrast, $\ALGORITHM$ is provably convergent, and furthermore is able to exploit asynchrony effectively to achieve speedup.


In addition, we also evaluate $\ALGORITHM$ under a wide range of communication delays. 
Due to space limitation, we only show performance in terms of gradient norm in Fig.~\ref{fig:simulation_vary_delay}.
We note that $\ALGORITHM$ converges in all cases, demonstrating its resilience against various delays in practice. 
Furthermore, as delay decreases, convergence becomes faster as robots have access to more up-to-date information from each other.

\subsection{Evaluation on benchmark PGO datasets}
\label{subsec:benchmark_results}


\begin{table*}[]
\setlength{\tabcolsep}{4pt}
\renewcommand{\arraystretch}{1.5}
\centering
\blue{
\begin{tabular}{ccc|cccc|cccc}
\hline\hline
\multirow{2}{*}{Dataset} & \multirow{2}{*}{\# Poses} & \multirow{2}{*}{\# Edges} & \multicolumn{4}{c|}{Cost value $f$} & \multicolumn{4}{c}{Additional statistics for ASAPP} \\ \cline{4-11} 
 &  &  & Opt.~\cite{Briales2017RAL} & Initial & DGS \cite{Choudhary2017Distributed} & ASAPP [ours] & Grad. Norm & Rot. RMSE {[}chordal{]} & Trans. RMSE {[}m{]} & Stepsize \\ \hline\hline
CSAIL (2D) & 1045 & 1171 & 31.47 & 36.10 & 31.54 & \textbf{31.51} & 0.36 & 0.0017 & 0.03 & 1.0 \\ \hline
Intel (2D) & 1228 & 1483 & 393.7 & 1205.9 & 394.6 & \textbf{393.7} & 0.002 & $2 \times 10^{-6}$ & $7 \times 10^{-5}$ & 1.0 \\ \hline
Manhattan (2D) & 3500 & 5453 & 193.9 & 1187.2 & 242.7 & \textbf{227.7} & 5.42 &0.07 & 1.18 & 0.03 \\ \hline
Garage (3D) & 1661 & 6275 & 1.267 & 4.477 & \textbf{1.277} & 1.309 & 0.04 & 0.014 & 0.41 & 0.05 \\ \hline
Sphere (3D) & 2500 & 4949 & 1687.0 & 3075.7 & 1743.1 & \textbf{1711.7} & 2.73 & 0.02 & 0.50 & 0.23 \\ \hline
Torus (3D) & 5000 & 9048 & 24227 & 25812 & 24305 & \textbf{24240} & 8.38 & 0.017 & 0.07 & 1.0 \\ \hline
Cubicle (3D) & 5750 & 16869 & 717.1 & 916.2 & \textbf{720.8} & 734.5 & 12.67 & 0.030 & 0.08 & 0.065 \\ \hline\hline
\end{tabular}
}
\caption{
    \small \blue{
	Evaluation on benchmark PGO datasets. 
	Each dataset is divided into trajectories of $5$ robots.
    We run $\ALGORITHM$ for $60$~s under a fixed communication delay of $0.1$~s.
    For reference, we also run DGS \cite{Choudhary2017Distributed} for $600$ synchronous iterations. 
    We compare the final cost values of the two approaches, and highlight the better solution in bold.
    For $\ALGORITHM$, we also report the used stepsize, achieved gradient norm, and rotation and translation root mean squared errors (RMSE) with respect to the global minimizer, computed by Cartan-Sync \cite{Briales2017RAL}.}}
\label{tab:slam_benchmarks}
\end{table*}

\blue{
We evaluate $\ALGORITHM$ on benchmark PGO datasets and compare its performance with Distributed Gauss-Seidel (DGS) \cite{Choudhary2017Distributed}, a state-of-the-art synchronous approach for distributed PGO used in recent multi-robot SLAM systems \cite{Cieslewski2018ICRA,Lajoie2019RAL}.
Each dataset is divided into $5$ segments simulating a collaborative SLAM mission with $5$ robots.

Following Choudhary et al.~\cite{Choudhary2017Distributed}, we initialize rotation estimates via distributed chordal initialization. 
To initialize the translations, we fix the rotation estimates and use distributed Gauss-Seidel to solve the reduced linear system over translations.
To test on scenarios where accurate initializations are not available, we restrict the number of Gauss-Seidel iterations to 50 for both rotation and translation initialization. 

Starting from the initial estimate, we run $\ALGORITHM$ with preconditioning for $60$s, assuming for simplicity a fixed delay of $0.1$s. 
Accordingly, we run DGS \cite{Choudhary2017Distributed} on the problem (linearized at the initial estimate) for $60/0.1=600$ synchronous iterations.
This setup favors DGS inherently, since each DGS iteration requires robots to communicate multiple times and update according to a specific order, which is likely to increase execution time in reality.

Table~\ref{tab:slam_benchmarks} compares the final cost values achieved by the two approaches,
where for $\ALGORITHM$ we first round the solutions to $\SE(d)$.
For $\ALGORITHM$, we also report the used stepsize, final gradient norm, and estimation errors with respect to the global minimizer computed by Cartan-Sync \cite{Briales2017RAL}.
As our results show, $\ALGORITHM$ often compares favorably against DGS, especially when the quality of initialization is poor. 
This is an important advantage, as good initialization schemes (such as distributed chordal initialization developed in \cite{Choudhary2017Distributed}) are usually iterative and thus expensive in terms of communication.
Furthermore, the $\ALGORITHM$ solution is close to the global minimizer, except on the \texttt{Manhattan} dataset where the rotation and translation errors are relatively high. 

We conclude this section by observing that on certain large datasets, convergence of $\ALGORITHM$ is slow as the iterate approaches a critical point. 
This is a consequence of the sublinear convergence rate, and in our case convergence is further impacted by the presence of communication delay. 
To address this, future work could consider accelerated methods to achieve higher precision.
To this end, the recent paper \cite{Fan2019ISRR} that generalizes Nesterov's accleration to PGO provides a promising direction.

}



\section{Conclusion}
We presented $\ALGORITHM$, the first \emph{asynchronous} and \emph{provably delay-tolerant} algorithm to solve distributed pose graph optimization and its rank-restricted semidefinite relaxations.
$\ALGORITHM$ enables each robot to run its local optimization process at a high rate, without waiting for updates from its peers over the network.
Assuming a worst-case bound on the communication delay, we established the global first-order convergence of $\ALGORITHM$,
and showed the existence of a convergent stepsize whose value depends on the worst-case delay and inherent problem sparsity.
When there is no delay, we further showed that this stepsize matches exactly with the corresponding constant in synchronous algorithms. 
Numerical evaluations on both simulation and real-world datasets confirm the advantages of $\ALGORITHM$ in reducing overall execution time, 
and demonstrate its resilience against a wide range of communication delay.

\blue{Our theoretical study in Section~\ref{sec:convergence} assumes a worst-case bounded delay.
Future work could consider less conservative strategies.
For example, the recent paper \cite{Tao2017NIPS} establishes convergence of asynchronous coordinate descent under \emph{unbounded} delay, where it is only assumed that the tail distribution of delay decays sufficiently fast. 
}
Another open question is conditions under which stronger performance guarantees may hold for first-order methods, e.g., second-order optimality.
Recent works have shown promising results towards this new direction \cite{Jin2017ICML,Criscitiello2019NIPS}.


\bibliographystyle{ieeetr}
\bibliography{paper}


\clearpage
\onecolumn
\appendix

\subsection{Proof of Lemma~\ref{lem:descent_lemma}}


\begin{proof}
Suppose that at iteration $k$, robot $\ik$ is selected for update.
Recall that $x^k = \begin{bmatrix}
x^k_1 & x^k_2 & \hdots & x^k_n
\end{bmatrix} \in \Mcal$, where $\Mcal$ is the product manifold $\Mcal = \Mcal_1 \times \Mcal_2 \times \hdots \Mcal_n$ (Section~\ref{subsec:global_view}). 
For all $k \in \mathbb{N}$, define the aggregate tangent vector $\eta^k \in T_{x^k} \Mcal$ as, 
\begin{equation}
    \eta^k_{i} \triangleq \begin{cases}
        \eta^k_{\ik}, & \text{if } i = \ik, \\
        0,  &\text{otherwise.}
    \end{cases}
\end{equation}
By definition of $\eta^k$, the update step in Algorithm~\ref{alg:global} (line~\ref{line:update_a}-\ref{line:update_b}) is equivalent to,
\begin{equation}
    x^{k+1} = \Retr_{x^k}(-\gamma \eta^k).
\end{equation}
By Lemma~\ref{lem:lipschitz_gradient_pullback}, $f$ has Lipschitz-type gradient for pullbacks. Therefore,
\begin{equation}
\begin{aligned}
    f(x^{k+1}) - f(x^k) &\leq -\gamma \langle \rgrad f(x^k), \eta^k \rangle + \frac{L\gamma^2}{2} \norm{\eta^k}^2 
    &= -\gamma \langle \rgrad_{\ik} f(x^k), \eta^k_{\ik} \rangle + \frac{L\gamma^2}{2} \norm{\eta^k_{\ik}}^2.
    \label{eq:descent_init_upper_bound}
\end{aligned}
\end{equation}
Using the equality $\langle \eta_1, \eta_2 \rangle = \frac{1}{2}[ \norm{\eta_1}^2 + \norm{\eta_2}^2 - \norm{\eta_1-\eta_2}^2]$, we can convert the inner product that appears on the right hand side of \eqref{eq:descent_init_upper_bound} into,
\begin{align}
    \langle \rgrad_{\ik} f(x^k), \eta^k_{\ik} \rangle &= 
    \frac{1}{2} \bigg [
    \norm{\rgrad_{\ik} f(x^k)}^2 +
    \norm{\eta^k_{\ik}}^2 -
    \norm{\rgrad_{\ik} f(x^k) - \eta^k_{\ik}}^2
    \bigg ].
    \label{eq:descent_inner_product}
\end{align}
Substitute \eqref{eq:descent_inner_product} into \eqref{eq:descent_init_upper_bound}. After collecting relevant terms, we have,
\begin{equation}
\begin{aligned}
    f(x^{k+1}) - f(x^k) &\leq -\gamma \langle \rgrad_{\ik} f(x^k), \eta^k_{\ik} \rangle + \frac{L\gamma^2}{2} \norm{\eta^k_{\ik}}^2 \\
    &= -\frac{\gamma}{2} \bigg [
    \norm{\rgrad_{\ik} f(x^k)}^2 +
    \norm{\eta^k_{\ik}}^2 -
    \norm{\rgrad_{\ik} f(x^k) - \eta^k_{\ik}}^2
    \bigg ]
    + 
    \frac{L\gamma^2}{2} \norm{\eta^k_{\ik}}^2
    \\
    &= -\frac{\gamma}{2} \norm{\rgrad_{\ik} f(x^k)}^2
    - \frac{\gamma-L\gamma^2}{2} \norm{\eta^k_{\ik}}^2
    + \frac{\gamma}{2} \norm{\rgrad_{\ik} f(x^k) - \eta^k_{\ik}}^2. 
    \label{eq:descent_upper_bound_2}
\end{aligned}
\end{equation}
We proceed to bound the last term on the right hand side of \eqref{eq:descent_upper_bound_2}. 
Recall from \eqref{eq:local_cost} and \eqref{eq:local_descent} that $\eta_{\ik}^k$ is formed with stale gradients,
\begin{align}
    \eta_{\ik}^k = \rgrad_{\ik} h_{\ik}(x^k_{\ik}) + \sum_{\jk \in \Nbr(\ik)}  \rgrad_{\ik} f_{\ek}(x^k_{\ik}, x^{k - B(\jk)}_{\jk}),
    \label{eq:outdated_gradient}
\end{align}
where we abbreviate the notation by defining $\ek \triangleq (\ik, \jk) \in E_{\Rcal}$.
In contrast, the Riemannian gradient $\rgrad_{\ik} f(x^k)$ is by definition formed using up-to-date variables,
\begin{align}
    \rgrad_{\ik} f(x^k) = \rgrad_{\ik} h_{\ik}(x^k_{\ik}) + \sum_{\jk \in \Nbr(\ik)}  \rgrad_{\ik} f_{\ek}(x^k_{\ik}, x^{k}_{\jk}). 
    \label{eq:up_to_date_gradient}
\end{align}
Note that the only difference between \eqref{eq:outdated_gradient} and \eqref{eq:up_to_date_gradient} is that 
delayed information is used in \eqref{eq:outdated_gradient}.
In order to form the last term on the right hand side of \eqref{eq:descent_upper_bound_2}, 
we subtract \eqref{eq:outdated_gradient} from \eqref{eq:up_to_date_gradient} and compute the norm distance.
Subsequently, we use the triangle inequality to obtain an upper bound on this norm distance,
\begin{equation}
    \begin{aligned}
    \norm{\rgrad_{\ik} f(x^k) - \eta^k_{\ik}} &= \norm{\sum_{\ek = (\ik,\jk) \in E_{\Rcal}} \bigg [ \rgrad_{\ik} f_{\ek}(x^k_{\ik}, x^{k}_{\jk}) - \rgrad_{\ik} f_{\ek}(x^k_{\ik}, x^{k - B(\jk)}_{\jk})  \bigg] } \\
    &\leq
    \sum_{\ek = (\ik,\jk) \in E_{\Rcal}} 
    \underbrace{
    \norm{\rgrad_{\ik} f_{\ek}(x^k_{\ik}, x^{k}_{\jk}) - \rgrad_{\ik} f_{\ek}(x^k_{\ik}, x^{k - B(\jk)}_{\jk})}}_{\epsilon(\ik,\jk)}.
    \label{eq:gradnorm_bound}
    \end{aligned}
\end{equation}
Next, we proceed to bound each $\epsilon(\ik,\jk)$ term. 
To do so, we use the fact that for a real-valued function $f$ defined over a matrix submanifold $\Mcal \subseteq \Ecal$, 
its Riemannian gradient is obtained as the orthogonal projection of the Euclidean gradient onto the tangent space (see \cite[Section.~3.6.1]{absil2009optimization}),
\begin{equation}
    \rgrad f(x) = \proj_{T_x \Mcal} \nabla f(x).
    \label{eq:egrad_to_rgrad}
\end{equation}
Substituting \eqref{eq:egrad_to_rgrad} into the right hand side of \eqref{eq:gradnorm_bound}, it holds that,
\begin{equation}
\begin{aligned}
    \epsilon(\ik,\jk) 
    &= 
    \norm{\rgrad_{\ik} f_{\ek}(x^k_{\ik}, x^{k}_{\jk}) - \rgrad_{\ik} f_{\ek}(x^k_{\ik}, x^{k - B(\jk)}_{\jk})} \\
    &=
    \norm{\proj_{T_{x_{\ik}^k}} \nabla_{\ik} f_{\ek}(x^k_{\ik}, x^{k}_{\jk}) - 
    \proj_{T_{x_{\ik}^k}} \nabla_{\ik} f_{\ek}(x^k_{\ik}, x^{k - B(\jk)}_{\jk})}.
    \label{eq:epsilon_bound_1}
\end{aligned}
\end{equation}
Furthermore, since the tangent space is identified as a linear subspace of the ambient space $\Ecal$ \cite[Section~3.5.7]{absil2009optimization}, 
the orthogonal projection operation is a \emph{non-expansive} mapping, i.e., 
\begin{equation}
\begin{aligned}
\epsilon(\ik,\jk) 
&= \norm{\proj_{T_{x_{\ik}^k}} \nabla_{\ik} f_{\ek}(x^k_{\ik}, x^{k}_{\jk}) - 
\proj_{T_{x_{\ik}^k}} \nabla_{\ik} f_{\ek}(x^k_{\ik}, x^{k - B(\jk)}_{\jk})} \\
&\leq
\norm{\nabla_{\ik} f_{\ek}(x^k_{\ik}, x^{k}_{\jk}) - 
      \nabla_{\ik} f_{\ek}(x^k_{\ik}, x^{k - B(\jk)}_{\jk})}.
\end{aligned}
\end{equation}
Since the norm distance with respect to $\ik$ is no greater than the overall norm distance, 
we furthermore have, 
\begin{equation}
    \epsilon(\ik,\jk) 
      \leq 
      \norm{\nabla_{\ik} f_{\ek}(x^k_{\ik}, x^{k}_{\jk}) - 
      \nabla_{\ik} f_{\ek}(x^k_{\ik}, x^{k - B(\jk)}_{\jk})}
      \leq 
      \norm{\nabla f_{\ek}(x^k_{\ik}, x^{k}_{\jk}) - 
          \nabla f_{\ek}(x^k_{\ik}, x^{k - B(\jk)}_{\jk})}.
\end{equation}
In \eqref{eq:P}, the Euclidean gradient of each cost function $f_{\ek}$ is Lipschitz continuous.
Furthermore, it is straightforward to show that the Lipschitz constant of $f_{\ek}$ is always less than or equal to the Lipschitz constant of the overall cost function $f$. 
Denote the latter as $C > 0$. 
By definition, we thus have, 
\begin{equation}
    \epsilon(\ik,\jk) \leq C\norm{x^{k}_{\jk} - x^{k - B(\jk)}_{\jk}}.
    \label{eq:epsilon_bound_2_1}
\end{equation}
In addition, in \cite{Tian2019RBCD} we have shown that the Riemannian version of the Lipschitz constant $L$ that appears in Lemma~\ref{lem:lipschitz_gradient_pullback} is always greater than or equal to the Euclidean Lipschitz constant $C$ (see Lemma 5 in \cite{Tian2019RBCD}). 
Thus,
\begin{equation}
    \epsilon(\ik,\jk) \leq L \norm{x^{k}_{\jk} - x^{k - B(\jk)}_{\jk}}.
    \label{eq:epsilon_bound_2}
\end{equation}
We proceed by bounding the norm on the right hand side of \eqref{eq:epsilon_bound_2}. 
Writing the subtraction as a telescoping sum and invoking triangle inequality, we first obtain,
\begin{equation}
\begin{aligned}
\norm{x^{k}_{\jk} - x^{k - B(\jk)}_{\jk}} = 
    \norm{\sum_{k' = k - B(\jk)}^{k-1} \bigg( 
    x_{\jk}^{k'+1} - x_{\jk}^{k'}
    \bigg)}
    \leq 
    \sum_{k' = k - B(\jk)}^{k-1} \norm{x_{\jk}^{k'+1} - x_{\jk}^{k'}}.
    \label{eq:epsilon_bound_3}
\end{aligned}
\end{equation}
Recall that for all $\jk$ and iterations $k'$, the next iterate $x_{\jk}^{k'+1}$ is obtained from $x_{\jk}^{k'}$ via the following update, 
\begin{equation}
    x_{\jk}^{k'+1} = \Retr_{x_{\jk}^{k'}} (-\gamma \eta_{\jk}^{k'}).
\end{equation}
Furthermore, from Lemma~5 in \cite{Tian2019RBCD}, we know that for each manifold $\Mcal_i$, there exists a corresponding constant $\alpha_i > 0$ such that the Euclidean distance from the initial point to the new point after retraction is always bounded by the following quantity, 
\begin{equation}
    \norm{\Retr_{x_{i}}(\eta_{i}) - x_{i}} \leq \alpha_{i} \norm{\eta_{i}}, \; \forall x \in \Mcal, \; \forall \eta_i \in T_{x_i} \Mcal.
    \label{eq:retraction_bound}
\end{equation}
Equation \eqref{eq:retraction_bound} thus provides a way to bound the term on the right hand side of \eqref{eq:epsilon_bound_3},
\begin{equation}
    \norm{x^{k}_{\jk} - x^{k - B(\jk)}_{\jk}} 
    \leq 
    \sum_{k' = k - B(\jk)}^{k-1}
    \norm{\Retr_{x_{\jk}^{k'}} (-\gamma \eta_{\jk}^{k'}) - x_{\jk}^{k'}}
    \leq
    \sum_{k' = k - B(\jk)}^{k-1}
    \alpha_{\jk} \norm{\gamma \eta_{\jk}^{k'}}.
\end{equation}
To remove the dependency on $\alpha_{\jk}$, let $\alpha \triangleq \max_{i \in \Rcal} \alpha_i$. We thus have,
\begin{equation}
    \norm{x^{k}_{\jk} - x^{k - B(\jk)}_{\jk}} \leq
    \alpha \gamma \sum_{k' = k - B(\jk)}^{k-1} \norm{\eta_{\jk}^{k'}}.
\end{equation}
We can further more use the bounded delay assumption (Assumption~\ref{as:bounded_delay}) to replace $B(\jk)$ with $B$,
\begin{align}
    \norm{x^{k}_{\jk} - x^{k - B(\jk)}_{\jk}} 
    \leq 
    \alpha \gamma \sum_{k' = k - B(\jk)}^{k-1} \norm{\eta_{\jk}^{k'}}
    \leq 
    \alpha  \gamma \sum_{k' = k - B}^{k-1} \norm{\eta_{\jk}^{k'}}.
    \label{eq:epsilon_bound_4}
\end{align}
Substituting \eqref{eq:epsilon_bound_4} into \eqref{eq:epsilon_bound_2}, we have,
\begin{align}
    \epsilon(\ik, \jk) \leq \alpha \gamma L \sum_{k' = k - B}^{k-1} \norm{\eta_{\jk}^{k'}}.
    \label{eq:epsilon_bound_5}
\end{align}
Substituting \eqref{eq:epsilon_bound_5} into \eqref{eq:gradnorm_bound}, we then have,
\begin{align}
    \norm{\rgrad_{\ik} f(x^k) - \eta_{\ik}^k} 
    \leq \sum_{\jk \in \Nbr(\ik)} \epsilon(\ik,\jk)
    \leq
    \alpha \gamma L \sum_{\jk \in \Nbr(\ik)} \sum_{k' = k - B}^{k-1} \norm{\eta_{\jk}^{k'}}.
    \label{eq:gradnorm_bound_2}
\end{align}
Squaring both sides of \eqref{eq:gradnorm_bound_2}, we obtain, 
\begin{equation}
    \norm{\rgrad_{\ik} f(x^k) - \eta_{\ik}^k}^2 
    \leq
    \bigg ( \alpha \gamma L \sum_{\jk \in \Nbr(\ik)} \sum_{k' = k - B}^{k-1} \norm{\eta_{\jk}^{k'}} \bigg)^2
    = \alpha^2 \gamma^2 L^2  \bigg (\sum_{\jk \in \Nbr(\ik)} \sum_{k' = k - B}^{k-1} \norm{\eta_{\jk}^{k'}} \bigg)^2.
    \label{eq:gradnorm_bound_22}
\end{equation}
Recall that the sum of squares inequality states that $(\sum_{i=1}^n a_i)^2 \leq n \sum_{i=1}^n a_i^2$. This gives an
upper bound on the squared term in \eqref{eq:gradnorm_bound_22},
\begin{equation}
\begin{aligned}
    \norm{\rgrad_{\ik} f(x^k) - \eta_{\ik}^k}^2 
    &\leq 
     \alpha^2 \gamma^2 L^2 B \Delta_{\ik} \sum_{\jk \in \Nbr(\ik)} \sum_{k' = k - B}^{k-1} \norm{\eta_{\jk}^{k'}}^2 
    &\leq 
     \alpha^2 \gamma^2 L^2 B \Delta \sum_{\jk \in \Nbr(\ik)} \sum_{k' = k - B}^{k-1} \norm{\eta_{\jk}^{k'}}^2,
    \label{eq:gradnorm_bound_3}
\end{aligned}
\end{equation}
where $\Delta_{\ik} \leq \Delta$ is robot $\ik$'s degree in the robot-level graph $G_{\Rcal}$.
Finally, substituting \eqref{eq:gradnorm_bound_3} in \eqref{eq:descent_upper_bound_2} concludes the proof,
\begin{align}
    f(x^{k+1}) - f(x^k) &\leq
    -\frac{\gamma}{2} \norm{\rgrad_{\ik} f(x^k)}^2
    - \frac{\gamma-L\gamma^2}{2} \norm{\eta^k_{\ik}}^2
    + \frac{\gamma}{2} \norm{\rgrad_{\ik} f(x^k) - \eta^k_{\ik}}^2  \\
    &\leq
    -\frac{\gamma}{2} \norm{\rgrad_{\ik} f(x^k)}^2
    - \frac{\gamma-L\gamma^2}{2} \norm{\eta^k_{\ik}}^2
    + \frac{ \alpha^2 \gamma^3 L^2 B \Delta}{2} 
    \sum_{\jk \in \Nbr(\ik)} \sum_{k' = k - B}^{k-1} \norm{\eta_{\jk}^{k'}}^2.
    \label{eq:descent_upper_bound_3}
\end{align}
\end{proof}


\subsection{Proof of Theorem~\ref{thm:exact_convergence}}

\begin{proof}
Since $f^\star$ is a global lower bound on $f$, we can obtain the following inequality,
\begin{equation}
    f^\star - f(x^0) \leq \mathbb{E}_{i_{0:K-1}} \bigg [f(x^K) \bigg] - f(x^0) 
    = \mathbb{E}_{i_{0:K-1}} \bigg [ 
    \sum_{k = 0}^{K-1} (f(x^{k+1}) - f(x^k))
    \bigg],
\end{equation}
where the right hand side rewrites the middle term as a telescoping sum. 
Using the linearity of expectation, we obtain,
\begin{equation}
     f^\star - f(x^0) \leq \mathbb{E}_{i_{0:K-1}} \bigg [ 
    \sum_{k = 0}^{K-1} (f(x^{k+1}) - f(x^k)) \bigg]
    = \sum_{k=0}^{K-1} 
    \mathbb{E}_{i_{0:k}} \bigg[f(x^{k+1}) - f(x^k)\bigg].
\end{equation}
For each expectation term, applying the law of total expectation yields,
\begin{equation}
    f^\star - f(x^0) \leq \sum_{k=0}^{K-1} 
    \mathbb{E}_{i_{0:k-1}} \bigg [ \mathbb{E}_{i_k|i_{0:k-1}}
    [f(x^{k+1}) - f(x^k)]
    \bigg].
    \label{eq:law_of_total_expectation}
\end{equation}
Next, recall that Lemma~\ref{lem:descent_lemma} already gives an upper bound on the innermost term of \eqref{eq:law_of_total_expectation}.
Substituting this upper bound into \eqref{eq:law_of_total_expectation} gives,
\begin{equation}
\begin{aligned}
    f^\star - f(x^0) \leq \sum_{k=0}^{K-1} 
    \mathbb{E}_{i_{0:k-1}} \bigg[
    \mathbb{E}_{i_k|i_{0:k-1}} \bigg [
    -\frac{\gamma}{2} \norm{\rgrad_{\ik} f(x^k)}^2
    - \frac{\gamma-L\gamma^2}{2} \norm{\eta^k_{\ik}}^2  
    + \frac{\Delta B \alpha^2 L^2 \gamma^3  }{2} 
    \sum_{\jk \in \Nbr(\ik)} \sum_{k' = k - B}^{k-1} \norm{\eta_{\jk}^{k'}}^2
    \bigg ]\bigg].
    \label{eq:convergence_thm_bound_1}
\end{aligned}
\end{equation}
Next, we simplify individual terms on the right hand side of \eqref{eq:convergence_thm_bound_1}.
We start with the first conditional expectation term. 
Using the definition of conditional expectation, 
\begin{equation}
    \mathbb{E}_{i_k|i_{0:k-1}} \bigg [-\frac{\gamma}{2} \norm{\rgrad_{\ik} f(x^k)}^2 \bigg ]
    = -\frac{\gamma}{2} \sum_{i = 1}^n P(\ik = i|i_{0:k-1}) \norm{\rgrad_{i} f(x^k)}^2.
\end{equation}
Recall that $\{\ik\}$ are i.i.d. random variables uniformly distributed over $1$ to $n$ (Section~\ref{subsec:global_view}).
Setting $P(\ik = i|i_{0:k-1}) = 1/n$ thus gives,
\begin{equation}
    \mathbb{E}_{i_k|i_{0:k-1}} \bigg [-\frac{\gamma}{2} \norm{\rgrad_{\ik} f(x^k)}^2 \bigg ]
    = -\frac{\gamma}{2} \sum_{i = 1}^n \frac{1}{n} \norm{\rgrad_{i} f(x^k)}^2
    = -\frac{\gamma}{2n} \norm{\rgrad f(x^k)}^2.
    \label{eq:conditional_expectation_1}
\end{equation}
Similarly, for the third conditional expectation in \eqref{eq:convergence_thm_bound_1}, we note that,
\begin{equation}
    \mathbb{E}_{i_k|i_{0:k-1}} \bigg [ 
    \sum_{\jk \in \Nbr(\ik)} \sum_{k' = k - B}^{k-1} \norm{\eta_{\jk}^{k'}}^2 \bigg ]
    =
    \sum_{i=1}^n \frac{1}{n}
    \sum_{j \in \Nbr(i)} \sum_{k' = k - B}^{k-1} \norm{\eta_{j}^{k'}}^2.
    \label{eq:conditional_expectation_2}
\end{equation}
In equation \eqref{eq:conditional_expectation_2}, exchange the order of summations and collect relevant terms,
\begin{equation}
    \sum_{i=1}^n \frac{1}{n}
    \sum_{j \in \Nbr(i)} \sum_{k' = k - B}^{k-1} \norm{\eta_{j}^{k'}}^2
    =
    \frac{1}{n} \sum_{k' = k - B}^{k-1} \sum_{i=1}^n  \sum_{j \in \Nbr(i)} \norm{\eta_{j}^{k'}}^2
    =
    \frac{2}{n} \sum_{k' = k - B}^{k-1} \norm{\eta^{k'}}^2.
    \label{eq:conditional_expectation_3}
\end{equation}
Using our simplified expressions for the first and third term on the right hand side of \eqref{eq:convergence_thm_bound_1}, we obtain,
\begin{align}
    f^\star - f(x^0) \leq \sum_{k=0}^{K-1} 
    \mathbb{E}_{i_{0:k-1}} \bigg[
    -\frac{\gamma}{2n} \norm{\rgrad f(x^k)}^2
    - \mathbb{E}_{i_k|i_{0:k-1}} \bigg[\frac{\gamma-L\gamma^2}{2} \norm{\eta^k_{\ik}}^2\bigg] 
    + \frac{\Delta B \alpha^2 L^2 \gamma^3}{n} 
    \sum_{k' = k - B}^{k-1} \norm{\eta^{k'}}^2
    \bigg].
    \label{eq:convergence_thm_bound_2}
\end{align}
Next, using the independence relations and the linearity of expectation, we obtain, 
\begin{align}
    f^\star - f(x^0) \leq \mathbb{E}_{i_{0:K-1}}
    \sum_{k=0}^{K-1} 
    \bigg[ 
    -\frac{\gamma}{2n} \norm{\rgrad f(x^k)}^2
    -\frac{\gamma-L\gamma^2}{2} \norm{\eta^k}^2
    +\frac{\Delta B \alpha^2 L^2 \gamma^3}{n} 
    \sum_{k' = k - B}^{k-1} \norm{\eta^{k'}}^2
    \bigg].
    \label{eq:convergence_thm_bound_3}
\end{align}
At this point, our bound still involves the squared norms of update vectors from earlier iterations (last term on the right hand side).
To simplify the bound further, note that,
\begin{align}
    \sum_{k=0}^{K-1} \sum_{k' = k - B}^{k-1} \norm{\eta^{k'}}^2 \leq B \sum_{k=0}^{K-1} \norm{\eta^{k}}^2.
\end{align}
Using the above inequality in \eqref{eq:convergence_thm_bound_3}, we obtain, 
\begin{align}
    f^\star - f(x^0) \leq \mathbb{E}_{i_{0:K-1}}
    \sum_{k=0}^{K-1} 
    \bigg[ 
    -\frac{\gamma}{2n} \norm{\rgrad f(x^k)}^2
    + \underbrace{(\frac{\Delta  \alpha^2 B^2 L^2 \gamma^3}{n} + \frac{L\gamma^2 - \gamma}{2})}_{A_1(\gamma)} \norm{\eta^{k}}^2
    \bigg].
    \label{eq:convergence_thm_bound_4}
\end{align}
We establish a sufficient condition on $\gamma$ such that $A_1(\gamma)$ as a whole is nonpositive.
Let us define $\rho \triangleq \Delta/n$. Consider the following factorization of $A_1(\gamma)$,
\begin{equation}
A_1(\gamma) = \frac{\gamma}{2} \underbrace{(2 \rho \alpha^2 B^2 L^2 \gamma^2 + L \gamma - 1)}_{A_2(\gamma)}.
\end{equation}
Note that $A_2(\gamma)$ is the same as \eqref{eq:stepsize_inequality} in Theorem~\ref{thm:exact_convergence}.
For the moment, suppose that we can find $\gamma > 0$ such that $A_2(\gamma) \leq 0$. 
This implies that $A_1(\gamma) \leq 0$, and we can thus drop this term on the right hand side of \eqref{eq:convergence_thm_bound_4},
\begin{equation}
\begin{aligned}
    f^\star - f(x^0) &\leq 
    -\frac{\gamma}{2n} 
    \sum_{k=0}^{K-1} \mathbb{E}_{i_{0:K-1}} \bigg[\norm{\rgrad f(x^k)}^2 \bigg].
\end{aligned}
\end{equation}
Replacing the expected value at each iteration with the minimum across all iterations, we have,
\begin{equation}
    f^\star - f(x^0) \leq -\frac{\gamma K}{2n} 
    \min_{k \in \{0,\hdots,K-1\}} \mathbb{E}_{i_{0:K-1}} \bigg[\norm{\rgrad f(x^k)}^2 \bigg].
    \label{eq:convergence_thm_bound_5}
\end{equation}
Finally, rearranging the last inequality yields, 
\begin{equation}
    \min_{k \in \{0,\hdots,K-1\}} \mathbb{E}_{i_{0:K-1}} \bigg[\norm{\rgrad f(x^k)}^2 \bigg] \leq 
    \frac{2n(f(x^0) - f^\star)}{\gamma K}.
    \label{eq:convergence_rate_appendix}
\end{equation}
Thus we have proved the main part of Theorem~\ref{thm:exact_convergence}. 
To prove the expression \eqref{eq:stepsize_bound}, note that if $B=0$ (i.e., there is no delay), the condition $A_2(\gamma) \leq 0$ entails $L\gamma \leq 1$. In this case, we can thus set the upper bound $\bar{\gamma}$ to $1/L$.
On the other hand, if $B > 0$, $A_2(\gamma)$ becomes a quadratic function of $\gamma$, and furthermore has the following positive root,
\begin{equation}
    \bar{\gamma} = \frac
    {\sqrt{1 + 8 \rho \alpha^2 B^2} - 1}
    {4 \rho \alpha^2 B^2 L} > 0.
\end{equation}
It is straightforward to verify that $A_2(\gamma) \leq 0$ for all $\gamma \in (0, \bar{\gamma}]$.
Therefore, we have proved that the condition $A_2(\gamma) \leq 0$ is ensured by the following choice of $\bar{\gamma}$, 
\begin{equation}
    \bar{\gamma} =
    \begin{cases}
    \frac
    {\sqrt{1 + 8 \rho \alpha^2 B^2 } - 1}
    {4 \rho \alpha^2 B^2 L}, & B > 0, \\
    1 / L, & B = 0.
    \end{cases}
\end{equation}
In particular, under this choice, $\ALGORITHM$ converges globally to first-order critical points, with an associated convergence rate given in \eqref{eq:convergence_rate_appendix}.

\end{proof}

\clearpage
\begin{figure*}[t]
	\centering
	\begin{subfigure}[t]{0.45\textwidth}
		\centering
		\includegraphics[width=\textwidth]{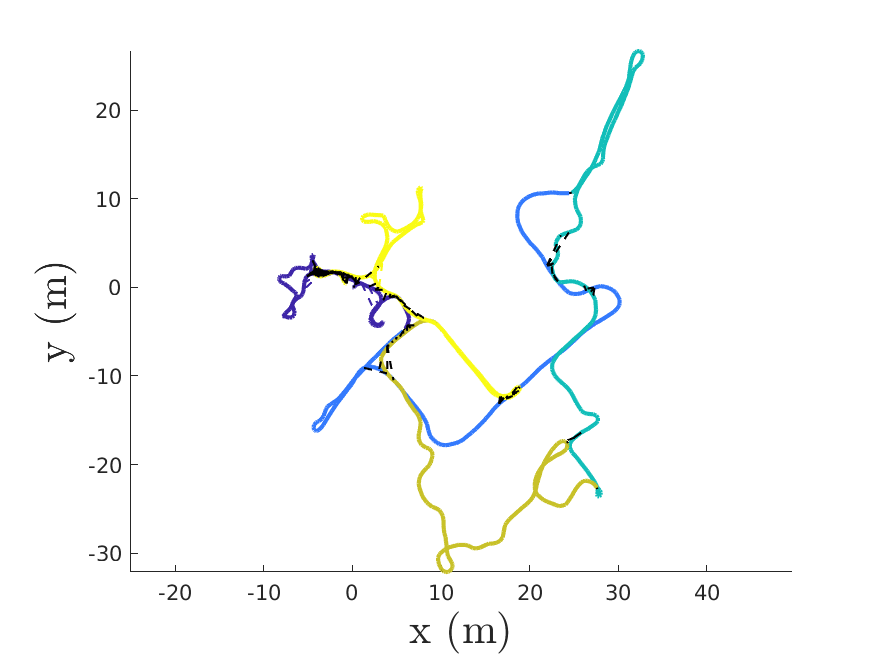}
		\caption{\small \texttt{CSAIL}}
		\label{fig:CSAIL}
	\end{subfigure}
	~
	\begin{subfigure}[t]{0.45\textwidth}
		\centering
		\includegraphics[width=\textwidth]{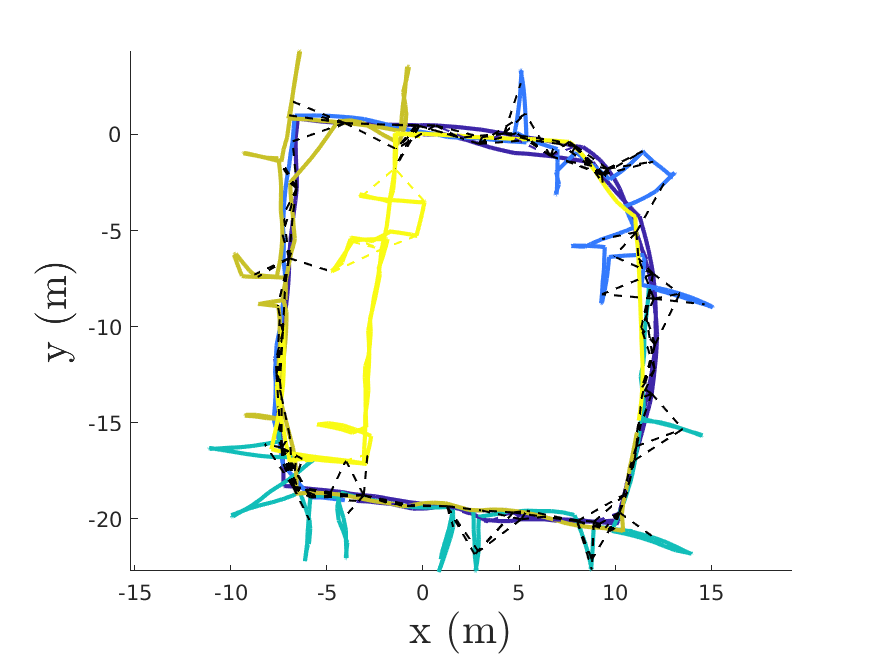}
		\caption{\small \texttt{Intel Research Lab}}
		\label{fig:intel}
	\end{subfigure}
	\\
	\begin{subfigure}[t]{0.45\textwidth}
		\centering
		\includegraphics[width=\textwidth]{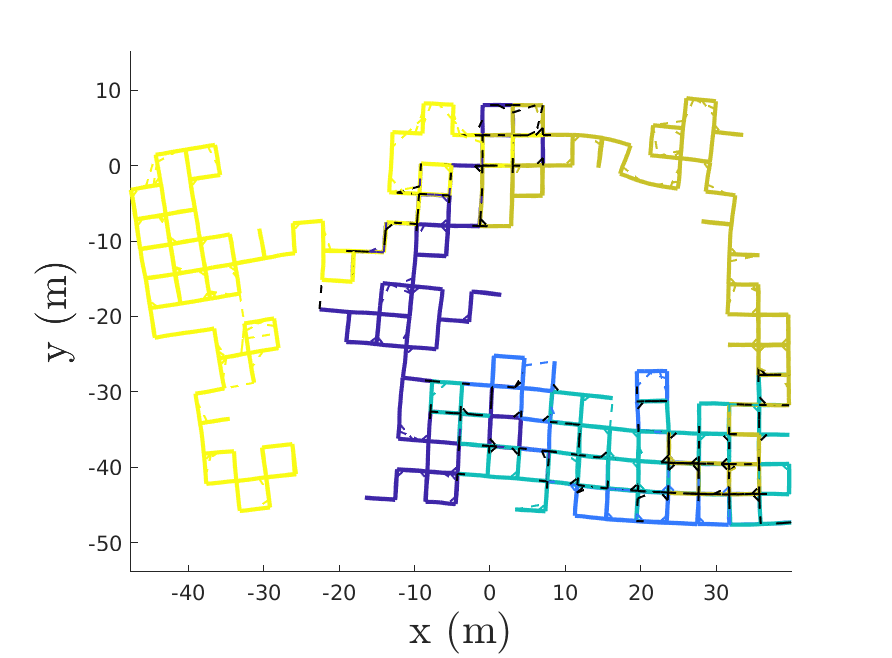}
		\caption{\small \texttt{Manhattan}}
		\label{fig:manhattan}
	\end{subfigure}
	~
	\begin{subfigure}[t]{0.45\textwidth}
		\centering
		\includegraphics[width=\textwidth]{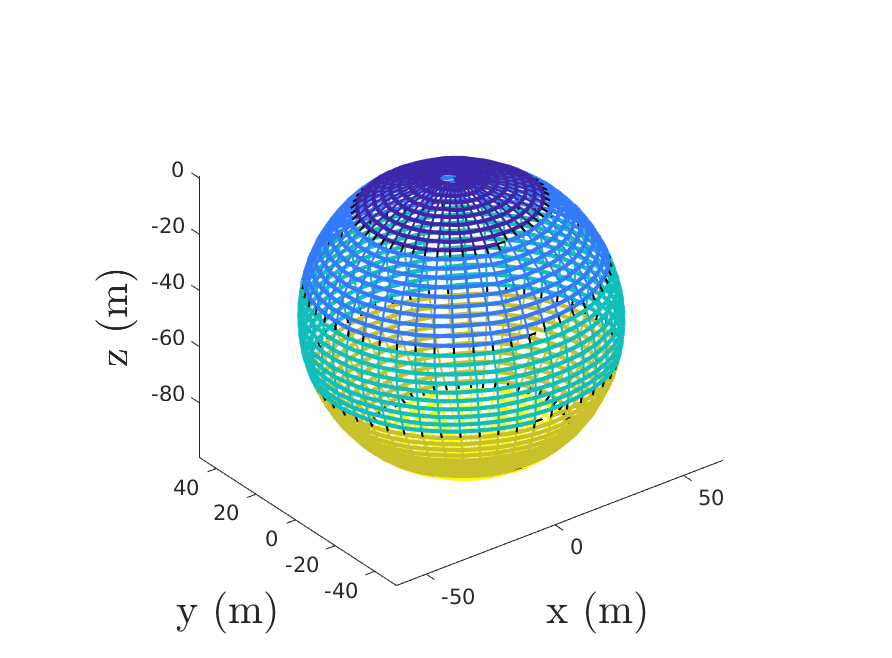}
		\caption{\small \texttt{Sphere}}
		\label{fig:sphere}
	\end{subfigure}
	\\
	\begin{subfigure}[t]{0.45\textwidth}
		\centering
		\includegraphics[trim=60 0 0 100, clip, width=\textwidth]{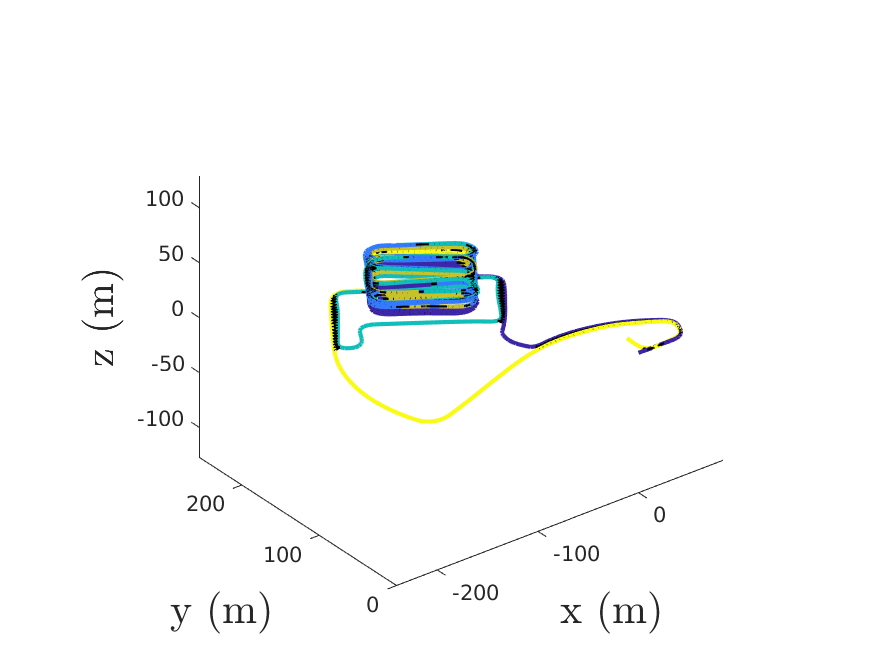}
		\caption{\small \texttt{Parking Garage}}
		\label{fig:parking_garage}
	\end{subfigure}
	~
	\begin{subfigure}[t]{0.45\textwidth}
		\centering
		\includegraphics[width=\textwidth]{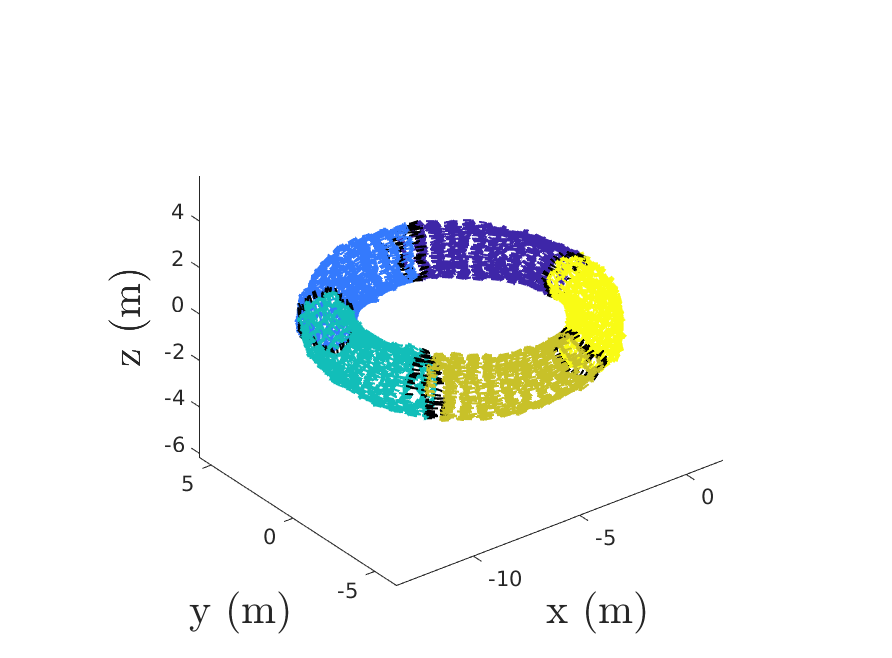}
		\caption{\small \texttt{Torus}}
		\label{fig:torus}
	\end{subfigure}
	\caption{\small 
	    Trajectory estimates returned by $\ALGORITHM$ on several benchmark datasets.
	    Each dataset contains trajectories of 5 robots (different colors). 
	    Inter-robot measurements (loop closures) are shown as black dashed lines.
	    (a) \texttt{CSAIL};
	    (b) \texttt{Intel Research Lab};
	    (c) \texttt{Manhattan};
	    (d) \texttt{Sphere};
	    (e) \texttt{Parking Garage};
	    (f) \texttt{Torus}.}
	\label{fig:datasets}
\end{figure*}

\end{document}